\documentclass[11pt,leqno]{article}
\usepackage{amsfonts}
\usepackage{amsmath}
\usepackage{amssymb}
\usepackage{theorem}
\usepackage{latexsym}
\usepackage{mathrsfs}

\newtheorem{theorem}{\bf Theorem}
\newtheorem{lemma}[theorem]{\bf Lemma}

\newtheorem{proposition}[theorem]{\bf Proposition}
\newtheorem{remark}[theorem]{\bf Remark}

\def\thebibliography#1{\section*{References\markboth
 {REFERENCES}{REFERENCES}}\list
 {[\arabic{enumi}]}{\settowidth\labelwidth{[#1]}\leftmargin\labelwidth
 \advance\leftmargin\labelsep
 \usecounter{enumi}}
 \def\newblock{\hskip .11em plus .33em minus .07em}
 \sloppy
 \sfcode`\.=1000\relax}

\begin{document}
%------------------------------------------------------------------------------
\vspace*{0ex}
\begin{center}
{\Large\bf
Solvability of the initial value problem \\[0.5ex]
to the Isobe--Kakinuma model for water waves 
}
\end{center}

\begin{center}
Ryo Nemoto and Tatsuo Iguchi
\end{center}

\begin{abstract}
We consider the initial value problem to the Isobe--Kakinuma model for water waves and the 
structure of the model. 
The Isobe--Kakinuma model is the Euler--Lagrange equations for an approximate Lagrangian 
which is derived from Luke's Lagrangian for water waves by approximating the velocity 
potential in the Lagrangian. 
The Isobe--Kakinuma model is a system of second order partial differential equations and 
is classified into a system of nonlinear dispersive equations. 
Since the hypersurface $t=0$ is characteristic for the Isobe--Kakinuma model, 
the initial data have to be restricted in an infinite dimensional manifold for 
the existence of the solution. 
Under this necessary condition and a sign condition, which corresponds to 
a generalized Rayleigh--Taylor sign condition for water waves, on the initial data, 
we show that the initial value problem is solvable locally in time in Sobolev spaces. 
We also discuss the linear dispersion relation to the model. 
\end{abstract}

%------------------------------------------------------------------------------
\section{Introduction}
\label{section:intro}
In this paper we are concerned with the solvability locally in time of the initial value problem 
to the Isobe--Kakinuma model for water waves. 
The water wave problem is mathematically formulated as a free boundary problem for 
an irrotational flow of an inviscid and incompressible fluid under the gravitational field. 
We consider the water filled in $(n+1)$-dimensional Euclidean space. 
Let $t$ be the time, $x=(x_1,\ldots,x_n)$ the horizontal spatial coordinates, 
and $z$ the vertical spatial coordinate. 
We assume that the water surface and the bottom are represented as $z=\eta(x,t)$ and $z=-h+b(x)$, 
respectively, where $\eta=\eta(x,t)$ is the surface elevation, $h$ is the mean depth, and $b=b(x)$ 
represents the bottom topography. 
J. C. Luke \cite{Luke1967} showed that the water wave problem has a variational structure, 
that is, he gave a Lagrangian in terms of the velocity potential $\Phi=\Phi(x,z,t)$ and the surface 
elevation $\eta$ in the form 
\begin{equation}\label{intro:Luke's Lagrangian}
\mathscr{L}_{\rm Luke}(\Phi,\eta) = \int_{-h+b(x)}^{\eta(x,t)}\biggl(\partial_t\Phi(x,z,t)
 +\frac12|\nabla_X\Phi(x,z,t)|^2+gz\biggr){\rm d}z
\end{equation}
and the action function 
\[
\mathscr{J}(\Phi,\eta)
= \int_{t_0}^{t_1}\!\!\!\int_{\Omega}\mathscr{L}_{\rm Luke}(\Phi,\eta){\rm d}x{\rm d}t, 
\]
where $\nabla_X=(\nabla,\partial_z)=(\partial_{x_1},\ldots,\partial_{x_n},\partial_z)$, 
$g$ is the gravitational constant, and $\Omega$ is an appropriate region in $\mathbf{R}^n$. 
He showed that the corresponding Euler--Lagrange equation is exactly the basic 
equations for water waves. 
M. Isobe \cite{Isobe1994, Isobe1994-2} and T. Kakinuma \cite{Kakinuma2000, Kakinuma2001, Kakinuma2003} 
approximated the velocity potential $\Phi$ in Luke's Lagrangian by 
\[
\Phi^{\mbox{\rm\tiny app}}(x,z,t) = \sum_{i=0}^N\Psi_i(z;b)\phi_i(x,t),
\]
where $\{\Psi_i\}$ is an appropriate function system in the vertical coordinate $z$ and may depend on 
the bottom topography $b$ and $(\phi_0,\phi_1,\ldots,\phi_N)$ are unknown variables, 
and derived an approximate Lagrangian 
$\mathscr{L}^{\mbox{\rm\tiny app}}(\phi_0,\phi_1,\ldots,\phi_N,\eta)
=\mathscr{L}_{\rm Luke}(\Phi^{\mbox{\rm\tiny app}},\eta)$. 
The Isobe--Kakinuma model is the corresponding Euler--Lagrange equation for the 
approximated Lagrangian. 
Different choices of the function system $\{\Psi_i\}$ give different Isobe--Kakinuma models. 
In this paper we adopt the approximation 
\begin{equation}\label{intro:approximation}
\Phi^{\mbox{\rm\tiny app}}(x,z,t) = \sum_{i=0}^N(z+h-b(x))^{p_i}\phi_i(x,t),
\end{equation}
where $p_0,p_1,\ldots,p_N$ are nonnegative integers satisfying $0=p_0<p_1<\cdots<p_N$. 
As we will see later, a natural choice of these exponents is given by 
$p_i=2i$ in the case of the flat bottom, that is the case $b(x)\equiv0$, 
and $p_i=i$ in the case of variable bottom topographies. 
In order to treat the both cases at the same time, we consider such a general case. 
Plugging this into Luke's Lagrangian \eqref{intro:Luke's Lagrangian} we obtain an 
approximate Lagrangian 
\begin{align*}
& \mathscr{L}^{\mbox{\rm\tiny app}}(\phi_0,\phi_1,\ldots,\phi_N,\eta) \\
&= \sum_{i=0}^N\frac{1}{p_i+1}H^{p_i+1}\partial_t\phi_i \\
&\quad\;
 +\frac12\sum_{i,j=0}^N\biggl(
  \frac{1}{p_i+p_j+1}H^{p_i+p_j+1}\nabla\phi_i\cdot\nabla\phi_j
  -\frac{2p_i}{p_i+p_j}H^{p_i+p_j}\phi_i\nabla b\cdot\nabla\phi_j \\
&\makebox[6em]{}
 +\frac{p_ip_j}{p_i+p_j-1}H^{p_i+p_j-1}(1+|\nabla b|^2)\phi_i\phi_j\biggr) \\
&\quad\;
 +\frac12g(\eta^2-(-h+b)^2),
\end{align*}
where $H=H(x,t)$ is the depth of the water and is given by $H(x,t)=h+\eta(x,t)-b(x)$. 
Here and in what follows we use the notational convention $0/0=0$. 
Then, the corresponding Euler--Lagrange equation has the form 
\begin{equation}\label{intro:IK model}
\left\{
 \begin{array}{l}
  \displaystyle
  H^{p_i}\partial_t\eta+\sum_{j=0}^N\biggl\{\nabla\cdot\biggl(
   \frac{1}{p_i+p_j+1}H^{p_i+p_j+1}\nabla\phi_j
   -\frac{p_j}{p_i+p_j}H^{p_i+p_j}\phi_j\nabla b\biggr) \\
  \displaystyle\phantom{ H^{p_i}\partial_t\eta+\sum_{j=0}^N\biggl\{ }
   +\frac{p_i}{p_i+p_j}H^{p_i+p_j}\nabla b\cdot\nabla\phi_j
   -\frac{p_ip_j}{p_i+p_j-1}H^{p_i+p_j-1}(1+|\nabla b|^2)\phi_j\biggr\}=0 \\
  \makebox[27em]{}\mbox{for}\quad i=0,1,\ldots,N, \\
  \displaystyle
  \sum_{i=0}^NH^{p_i}\partial_t\phi_i+g\eta 
   +\frac12\sum_{i,j=0}^N\Bigl(H^{p_i+p_j}\nabla\phi_i\cdot\nabla\phi_j
   -2p_iH^{p_i+p_j-1}\phi_i\nabla b\cdot\nabla\phi_j \\
  \displaystyle\makebox[12em]{}
   +p_ip_jH^{p_i+p_j-2}(1+|\nabla b|^2)\phi_i\phi_j\Bigr)=0.
 \end{array}
\right.
\end{equation}
This is the Isobe--Kakinuma model that we are going to consider in this paper. 
We consider the initial value problem to this Isobe--Kakinuma model \eqref{intro:IK model} 
under the initial conditions 
\begin{equation}\label{intro:initial conditions}
(\eta,\phi_0,\ldots,\phi_N)=(\eta_{(0)},\phi_{0(0)},\ldots,\phi_{N(0)}) \quad\makebox[3em]{at} t=0.
\end{equation}
Unique solvability locally in time of the initial value problem 
\eqref{intro:IK model}--\eqref{intro:initial conditions} in the case $N=1$ and $p_1=2$ 
and fundamental properties of the model are presented in 
Y. Murakami and T. Iguchi \cite{MurakamiIguchi2015}. 
Therefore, this paper is a generalization of their results.

One of the interesting features of the model is its linear dispersion relation. 
In the following section we will consider the linearized equations of the model around the rest state in the case 
of the flat bottom and calculate the linear dispersion relation together with the phase speed 
$c_{IK}(\xi)$ of the plane wave solution related to the wave vector $\xi\in\mathbf{R}^n$. 
See \eqref{dispersion:phase speed}. 
It is well known that the phase speed $c_{WW}(\xi)$ of the plane wave solution to the linearized equations 
for water waves is given by 
\[
c_{WW}(\xi)=\pm\sqrt{gh\frac{\tanh(h|\xi|)}{h|\xi|}}.
\]
If we choose $p_i=2i$, then we can show that 
\begin{equation}\label{intro:Pade}
(c_{IK}(\xi))^2 = [2N/2N] \mbox{ Pad\'e approximant of } (c_{WW}(\xi))^2,
\end{equation}
which will be given in Theorem \ref{dispersion:theorem 1}. 
Concerning Pad\'e approximants we refer to G. A. Baker and P. Graves-Morris \cite{BakerGraves-Morris1996}. 
This relation \eqref{intro:Pade} implies that the Isobe--Kakinuma model gives a good approximation of the 
basic equations for water waves in the shallow water regime $h|\xi|\ll1$. 
In fact, T. Iguchi \cite{Iguchi2017} gave a mathematically rigorous justification of the Isobe--Kakinuma 
model in the case of the flat bottom with the choice $N=1$ and $p_1=2$. 
He showed that the Isobe--Kakinuma model gives a higher order shallow water approximation for water waves 
with an error of order $O(\delta^6)$, where $\delta$ is a small nondimensional parameter defined as the 
ratio of the mean depth $h$ to the typical wave length. 
We note that the Green--Naghdi equations are known as a higher order shallow water approximation for water waves 
with an error of order $\delta^4$. 
Therefore, the Isobe--Kakinuma model gives a better approximation than the Green--Naghdi equations in the 
shallow water regime. 
Concerning the shallow water approximations and the rigorous justifications of the Green--Naghdi equations 
we refer to T. Iguchi \cite{Iguchi2009, Iguchi2011}, B. Alvarez-Samaniego and D. Lannes 
\cite{Alvarez-SamaniegoLannes2008}, Y. A. Li \cite{Li2006}, H. Fujiwara and T. Iguchi \cite{FujiwaraIguchi2015}, 
and D. Lannes \cite{Lannes2013-2}. 
The relation \eqref{intro:Pade} anticipates that the Isobe--Kakinuma model \eqref{intro:IK model} in the 
case of the flat bottom with the choice $p_i=2i$ would give an approximation with an error of order 
$O(\delta^{4N+2})$. We postpone this subject in the future research. 
If we choose $p_i=i$, we do not have such a beautiful relation as \eqref{intro:Pade} any more. 
However, this choice of $p_i$ would be important in the case of the variable bottom topographies and 
we still have a good approximation, which will be stated in Theorem \ref{dispersion:theorem 2}.

The Isobe--Kakinuma model \eqref{intro:IK model} is written in the matrix form as 
\[
\left(
 \begin{array}{cccc}
  H^{p_0} & 0 & \cdots & 0 \\
  \vdots  & \vdots & & \vdots \\
  H^{p_N} & 0 & \cdots & 0 \\
  0 & H^{p_0} & \cdots & H^{p_N}
 \end{array}
\right)
\partial_t
\left(
 \begin{array}{c}
  \eta \\
  \phi_0 \\
  \vdots \\
  \phi_N
 \end{array}
\right)
+\{\mbox{spatial derivatives}\}=\mbox{\boldmath$0$}.
\]
Since the coefficient matrix always has the zero eigenvalue, the hypersurface $t=0$ in the space-time 
$\mathbf{R}^n\times\mathbf{R}$ is characteristic for the Isobe--Kakinuma model \eqref{intro:IK model}, 
so that the initial value problem \eqref{intro:IK model}--\eqref{intro:initial conditions} is not solvable in general. 
In fact, if the problem has a solution $(\eta,\phi_0,\ldots,\phi_N)$, then by eliminating the time derivative 
$\partial_t\eta$ from the equations we see that the solution has to satisfy the relation 
\begin{align}\label{intro:compatibility}
& H^{p_i}\sum_{j=0}^N\nabla\cdot\biggl(
   \frac{1}{p_j+1}H^{p_j+1}\nabla\phi_j
   -\frac{p_j}{p_j}H^{p_j}\phi_j\nabla b\biggr) \nonumber \\
&= \sum_{j=0}^N\biggl\{\nabla\cdot\biggl(
   \frac{1}{p_i+p_j+1}H^{p_i+p_j+1}\nabla\phi_j
   -\frac{p_j}{p_i+p_j}H^{p_i+p_j}\phi_j\nabla b\biggr)  \\
&\phantom{ =\sum_{j=0}^N\biggl\{ }
 \displaystyle
   +\frac{p_i}{p_i+p_j}H^{p_i+p_j}\nabla b\cdot\nabla\phi_j
   -\frac{p_ip_j}{p_i+p_j-1}H^{p_i+p_j-1}(1+|\nabla b|^2)\phi_j\biggr\} \nonumber 
\end{align}
for $i=1,\ldots,N$. 
Therefore, as a necessary condition the initial date $(\eta_{(0)},\phi_{0(0)},\ldots,\phi_{N(0)})$ and 
the bottom topography $b$ have to satisfy the relation \eqref{intro:compatibility} for the existence of 
the solution. 
As we will see in Proposition \ref{intro:prop}, 
such initial data are constructed from the three scalar functions $b$, $\eta_{(0)}$, 
and $\phi_{(0)}$, where $\phi_{(0)}$ is the initial data for the trace of the velocity potential on the water surface.

The water wave problem has a conserved energy 
\[
E_{WW}(t) = \frac{\rho}{2}\int_{\mathbf{R}^n}\biggl\{\int_{-h+b(x)}^{\eta(x,t)}\bigl|
 \nabla_X\Phi(x,z,t)\bigr|^2{\rm d}z 
 +g\bigl(\eta(x,t)\bigr)^2\biggr\}{\rm d}x, 
\]
where $\rho$ is the constant density of the water. 
The first term in the right-hand side is the kinetic energy and the second one is the potential energy 
due to the gravity. 
The Isobe--Kakinuma model \eqref{intro:IK model} has also a conserved energy 
\begin{equation}\label{intro:energy}
E(t) = \frac{\rho}{2}\int_{\mathbf{R}^n}\biggl\{\int_{-h+b(x)}^{\eta(x,t)}\bigl|
 \nabla_X\Phi^{\mbox{\rm\tiny app}}(x,z,t)\bigr|^2{\rm d}z 
 +g\bigl(\eta(x,t)\bigr)^2\biggr\}{\rm d}x, 
\end{equation}
where $\Phi^{\mbox{\rm\tiny app}}$ is given by \eqref{intro:approximation}. 
See \eqref{proof:energy} for more explicit expression of this energy function $E(t)$. 
Under a physically reasonable condition on the water surface $\eta$ and the bottom topography $b$ 
we have an equivalence 
\[
E(t) \simeq \int_{\mathbf{R}^n}\biggl\{
 |\nabla\phi_0(x,t)|^2+\sum_{i=1}^N\Bigl(|\nabla\phi_i(x,t)|^2+\bigl(\phi_i(x,t)\bigr)^2\Bigr)
 +\bigl(\eta(x,t)\bigr)^2 \biggr\}{\rm d}x.
\]
Therefore, it is natural to work in the function space $\eta,\nabla\phi_0\in C([0,T];H^m)$ and 
$\phi_1,\ldots,\phi_N\in C([0,T];H^{m+1})$, where $H^m=W^{m,2}(\mathbf{R}^n)$ is the standard $L^2$ 
Sobolev space of order $m$ on $\mathbf{R}^n$.

It is well known that the well-posedness of the initial value problem to the water wave problem may be 
broken unless a generalized Rayleigh--Taylor sign condition $-\frac{\partial P}{\partial N} \geq c_0>0$ 
on the water surface is satisfied, where $P$ is the pressure and $N$ is the unit outward normal 
on the water surface. 
For example, we refer to S. Wu \cite{Wu1997, Wu1999} and D. Lannes \cite{Lannes2005}. 
This sign condition is equivalent to $-\partial_z P\geq c_0>0$ because the pressure $P$ is equal to 
the constant atmospheric pressure $P_0$ on the water surface. 
By using Bernoulli's law 
\[
\partial_t\Phi+\frac12|\nabla_X\Phi|^2+\frac{1}{\rho}(P-P_0)+gz\equiv0,
\]
the sign condition can be written in term of our unknowns $(\eta,\phi_0,\ldots,\phi_N)$ and $b$ as 
$a(x,t)\geq c_0>0$, where 
\begin{align}\label{intro:a}
a &= g+\sum_{i=0}^Np_iH^{p_i-1}\partial_t\phi_i \\
&\quad\;
+\frac12\sum_{i,j=0}^N\Bigl\{
 (p_i+p_j)H^{p_i+p_j-1}\nabla\phi_i\cdot\nabla\phi_j
 -2p_i(p_i+p_j-1)H^{p_i+p_j-2}\phi_i\nabla b\cdot\nabla\phi_j \nonumber \\
&\makebox[6em]{}
+p_ip_j(p_i+p_j-2)H^{p_i+p_j-3}(1+|\nabla b|^2)\phi_i\phi_j\Bigr\}. \nonumber
\end{align}
In fact, we have $-\frac{1}{\rho}\partial_z P^{\mbox{\rm\tiny app}}=g+\partial_z\partial_t\Phi^{\mbox{\rm\tiny app}}+
\nabla_X\partial_z\Phi^{\mbox{\rm\tiny app}}\cdot\nabla_X\Phi^{\mbox{\rm\tiny app}}=a$ on $z=\eta(x,t)$. 
In this paper, we assume that this sign condition is satisfied at the initial time $t=0$.

Now, we are ready to give our main result in this paper.

\begin{theorem}\label{intro:theorem}
Let $g, h, c_0, M_0$ be positive constants and $m$ an integer such that $m>n/2+1$. 
There exists a time $T>0$ such that if the initial data $(\eta_{(0)},\phi_{0(0)},\ldots,\phi_{N(0)})$ and 
$b$ satisfy the relation \eqref{intro:compatibility} and 
\begin{equation}\label{intro:conditions}
\left\{
 \begin{array}{l}
  \|\eta_{(0)}\|_m+\|\nabla\phi_{0(0)}\|_m+\|(\phi_{1(0)},\ldots,\phi_{N(0)})\|_{m+1}
   + \|b\|_{W^{m+2,\infty}}\leq M_0, \\[0.5ex]
  h+\eta_{(0)}(x)-b(x)\geq c_0, \quad a(x,0)\geq c_0 
   \qquad\mbox{for}\quad x\in\mathbf{R}^n,
 \end{array}
\right.
\end{equation}
then the initial value problem \eqref{intro:IK model}--\eqref{intro:initial conditions} 
has a unique solution $(\eta,\phi_0,\ldots,\phi_N)$ satisfying
$$
\eta,\nabla\phi_0\in C([0,T];H^m), \quad
  \phi_1,\ldots,\phi_N\in C([0,T];H^{m+1}).
$$
Moreover, the energy function $E(t)$ defined by \eqref{intro:energy} is a conserved quantity. 
\end{theorem}

\begin{remark}\label{intro:remark}
{\rm
(1) \ 
If we impose an additional condition $\phi_{(0)}\in L^2(\mathbf{R}^n)$, then the solution satisfies 
an additional integrability $\phi_0\in C([0,T];H^{m+1})$. 

(2) \ 
In the sign condition $a(x,0)\geq c_0>0$ we use the quantities 
$\partial_t\phi_1(x,0),\ldots,\partial_t\phi_N(x,0)$ which should be written in terms of the initial data. 
Although the hypersurface $t=0$ is characteristic for the Isobe--Kakinuma model \eqref{intro:IK model}, 
we can express $\partial_t\phi_1(x,0),\ldots,\partial_t\phi_N(x,0)$ in terms of the initial data and $b$. 
For details, we refer to Remark \ref{construction:remark}. 

(3) \ 
Let $\phi$ be the trace of the velocity potential $\Phi$ on the water surface. 
In view of our approximation \eqref{intro:approximation}, $\phi$ should be related to our variables 
approximately by the formula 
\begin{equation}\label{intro:relation}
\phi=\sum_{i=0}^NH^{p_i}\phi_i.
\end{equation}
Given the initial data $\eta_{(0)}$ and $\phi_{(0)} (=\phi|_{t=0})$ and the bottom topography $b$, 
the necessary condition \eqref{intro:compatibility} and the relation \eqref{intro:relation} determine 
uniquely the initial data $\phi_{0(0)},\ldots,\phi_{N(0)}$. 
In fact, we have the following proposition. 
}
\end{remark}

\begin{proposition}\label{intro:prop}
Let $h, c_0, M_0$ be positive constants and $m$ an integer such that $m>n/2+1$. 
There exists a positive $C>0$ such that if the initial data $(\eta_{(0)},\phi_{(0)})$ and $b$ satisfy 
\[
\left\{
 \begin{array}{l}
  \|\eta_{(0)}\|_m+\|b\|_{W^{m,\infty}}\leq M_0, \quad \|\nabla\phi_{(0)}\|_{m-1}<\infty, \\[0.5ex]
  h+\eta_{(0)}(x)-b(x)\geq c_0 \qquad\mbox{for}\quad x\in\mathbf{R}^n,
 \end{array}
\right.
\]
then the necessary condition \eqref{intro:compatibility} and the relation \eqref{intro:relation} determine 
uniquely the initial data $\phi_{0(0)},\ldots,\phi_{N(0)}$, which satisfy 
\[
\|\nabla\phi_{0(0)}\|_{m-1}+\|(\phi_{1(0)},\ldots,\phi_{N(0)})\|_m \leq C\|\nabla\phi_{(0)}\|_{m-1}. 
\]
\end{proposition}

\medskip
The contents of this paper are as follows. 
In Section \ref{section:dispersion} we consider the linearized equations of the Isobe--Kakinuma model 
around the rest state in the case of the flat bottom and analyze the linear dispersion relation. 
Especially, we show the beautiful relation \eqref{intro:Pade}. 
In Section \ref{section:linear} we consider the linearized equations of the Isobe--Kakinuma model 
around an arbitrary flow, reveal a hidden symmetric structure of the model, and give an idea to obtain 
an energy estimate for the solution of the nonlinear equations. 
In Section \ref{section:construction} we transform the Isobe--Kakinuma model to a system of equations 
for which the hypersurface $t=0$ is noncharacteristic by using the necessary condition \eqref{intro:compatibility} 
and construct the solution of the initial value problem to the transformed system by using a standard 
parabolic regularization. 
In Section \ref{section:proof} we show that the solution constructed in Section 
\ref{section:construction} is the solution of the Isobe--Kakinuma model \eqref{intro:IK model} 
if the initial data satisfy the necessary condition \eqref{intro:compatibility}.

\bigskip
\noindent
{\bf Notation}. \ 
We denote by $W^{m,p}(\mathbf{R}^n)$ the $L^p$ Sobolev space of order $m$ on $\mathbf{R}^n$. 
The norms of the Lebesgue space $L^p(\mathbf{R}^n)$ and the Sobolev space $H^m=W^{m,2}(\mathbf{R}^n)$ 
are denoted by $|\cdot|_p$ and $\|\cdot\|_m$, respectively. 
The $L^2$-norm and the $L^2$-inner product are simply denoted by $\|\cdot\|$ and $(\cdot,\cdot)_{L^2}$, 
respectively. 
We put $\partial_t=\partial/\partial t$, $\partial_j=\partial/\partial x_j$, 
and $\partial_z=\partial/\partial z$. 
For a multi-index $\alpha=(\alpha_1,\ldots,\alpha_n)$ we put 
$\partial^{\alpha}=\partial_1^{\alpha_1}\cdots\partial_n^{\alpha_n}$. 
$[P,Q]=PQ-QP$ denotes the commutator. 
For a matrix $A$ we denote by $A^{\rm T}$ the transpose of $A$. 
For a vector $\mbox{\boldmath$\phi$}=(\phi_0,\phi_1,\ldots,\phi_N)^{\rm T}$ we denote the last $N$ 
components by $\mbox{\boldmath$\phi$}'=(\phi_1,\ldots,\phi_N)^{\rm T}$. 
We use the notational convention $0/0=0$.

%------------------------------------------------------------------------------
\section{Linear dispersion relation}
\label{section:dispersion}
\setcounter{equation}{0}
\setcounter{theorem}{0}
In this section we consider the linearized equations of the Isobe--Kakinuma model \eqref{intro:IK model} 
in the case of the flat bottom. 
The linearized equations have the form 
\begin{equation}\label{dispersion:IK model}
\left\{
 \begin{array}{l}
  \displaystyle
  \partial_t\eta+\sum_{j=0}^N\biggl(
   \frac{1}{p_i+p_j+1}h^{p_j+1}\Delta\phi_j
   -\frac{p_ip_j}{p_i+p_j-1}h^{p_j-1}\phi_j\biggr)=0 \\
  \makebox[18em]{}\mbox{for}\quad i=0,1,\ldots,N, \\
  \displaystyle
  \sum_{i=0}^Nh^{p_i}\partial_t\phi_i+g\eta=0.
 \end{array}
\right.
\end{equation}
Putting $\mbox{\boldmath$\psi$}=(h^{p_0}\phi_0,\ldots,h^{p_N}\phi_N)^{\rm T}$ we can rewrite the above 
equations as the following simple matrix form 
\[
\left(
 \begin{array}{cc}
  0 & h\mbox{\boldmath$1$}^{\rm T} \\
  -h\mbox{\boldmath$1$} & O
 \end{array}
\right)
\partial_t
\left(
 \begin{array}{c}
  \eta \\
  \mbox{\boldmath$\psi$}
 \end{array}
\right)
+
\left(
 \begin{array}{cc}
  gh & \mbox{\boldmath$0$}^{\rm T} \\
  \mbox{\boldmath$0$} & A(hD)
 \end{array}
\right)
\left(
 \begin{array}{c}
  \eta \\
  \mbox{\boldmath$\psi$}
 \end{array}
\right)
= \mbox{\boldmath$0$},
\]
where $\mbox{\boldmath$1$}=(1,\ldots,1)^{\rm T}$ and $A(hD)=-A_0h^2\Delta+A_1$. 
The $(N+1)\times(N+1)$ matrices $A_0$ and $A_1$ are given by 
\[
A_0=\biggl(\frac{1}{p_i+p_j+1}\biggr)_{0\leq i,j\leq N}, \qquad
A_1=\biggl(\frac{p_ip_j}{p_i+p_j-1}\biggr)_{0\leq i,j\leq N},
\]
where we used a rather nonstandard notation for matrices. 
Since it might not cause any confusion, we will continue to use this notation in the following. 
Therefore, the linear dispersion relation is given by 
\[
\det
\left(
 \begin{array}{cc}
  gh & \sqrt{-1}h\omega\mbox{\boldmath$1$}^{\rm T} \\
  -\sqrt{-1}h\omega\mbox{\boldmath$1$} & A(h\xi)
 \end{array}
\right)=0,
\]
where $\xi\in\mathbf{R}^n$ is the wave vector, $\omega\in\mathbf{C}$ is the angular frequency, 
and $A(h\xi)=(h|\xi|)^2A_0+A_1$. 
We can rewrite the above dispersion relation as 
\begin{equation}\label{dispersion:dispersion}
h^2\omega^2\det\tilde{A}(h\xi)-gh\det A(h\xi)=0.
\end{equation}
Throughout this section we use the notation 
\[
\tilde{A}=
\left(
 \begin{array}{cc}
  0 & \mbox{\boldmath$1$}^{\rm T} \\
  -\mbox{\boldmath$1$} & A
 \end{array}
\right)
\]
for a matrix $A$. 
Concerning the determinants in the dispersion relation \eqref{dispersion:dispersion} we have the 
following proposition.

\begin{proposition}\label{dispersion:prop 1}
\begin{enumerate}
\setlength{\itemsep}{-3pt}
\item
For any $\xi\in\mathbf{R}^n\setminus\{\mbox{\boldmath$0$}\}$, the symmetric matrix $A(h\xi)$ is positive. 
\item
There exists $c_0>0$ such that for any $\xi\in\mathbf{R}^n$ we have $\det\tilde{A}(h\xi)\geq c_0$. 
\item
$(h|\xi|)^{-2}\det A(h\xi)$ is a polynomial in $(h|\xi|)^2$ of degree $N$ and the coefficient of $(h|\xi|)^{2N}$ 
is $\det A_0$. 
\item
$\det\tilde{A}(h\xi)$ is a polynomial in $(h|\xi|)^2$ of degree $N$ and the coefficient of $(h|\xi|)^{2N}$ 
is $\det\tilde{A}_0$. 
\end{enumerate}
\end{proposition}

\noindent
{\bf Proof}. \ 
For any $\mbox{\boldmath$\psi$}=(\psi_0,\ldots,\psi_N)^{\rm T}\in\mathbf{R}^{N+1}$ we see that 
\begin{align}\label{dispersion:identity}
(A(h\xi)\mbox{\boldmath$\psi$})\cdot\mbox{\boldmath$\psi$}
&= \int_0^1\biggl\{(h|\xi|)^2\biggl(\sum_{i=0}^N\psi_iz^{p_i}\biggr)^2
 +\biggl(\sum_{i=1}^Np_i\psi_iz^{p_i-1}\biggr)^2\biggr\}{\rm d}z \\
&\simeq (h|\xi|)^2|\mbox{\boldmath$\psi$}|^2+|\mbox{\boldmath$\psi$}'|^2, \nonumber
\end{align}
where $\mbox{\boldmath$\psi$}'=(\psi_1,\ldots,\psi_N)^{\rm T}$. 
This shows a positivity of $A(h\xi)$ for $\xi\ne\mbox{\boldmath$0$}$. 

Let $\xi\ne\mbox{\boldmath$0$}$ and $C(h\xi)=(C_{ij}(h\xi))_{0\leq i,j\leq N}$ be the cofactor matrix of $A(h\xi)$. 
Since $A(h\xi)$ is symmetric, $C(h\xi)$ is also symmetric and we have $C(h\xi)=(\det A(h\xi))A(h\xi)^{-1}$. 
Then, expanding the first row of the matrix $\tilde{A}(h\xi)$ we see that
\[
\det\tilde{A}(h\xi)=\sum_{i,j=0}^NC_{ij}(h\xi)
=(C(h\xi)\mbox{\boldmath$1$})\cdot\mbox{\boldmath$1$}
=(\det A(h\xi))(A(h\xi)^{-1}\mbox{\boldmath$1$})\cdot\mbox{\boldmath$1$},
\]
which implies the positivity of $\det\tilde{A}(h\xi)$ due to the positivity of $A(h\xi)$. 
We proceed to consider the case $\xi=\mbox{\boldmath$0$}$. 
It holds that $\tilde{A}(\mbox{\boldmath$0$})=\tilde{A}_1$ and 
\[
A_1=
\left(
 \begin{array}{cc}
  0 & \mbox{\boldmath$0$}^{\rm T} \\
  \mbox{\boldmath$0$} & A_1'
 \end{array}
\right), \quad
A_1'=\biggl(\frac{p_ip_j}{p_i+p_j-1}\biggr)_{1\leq i,j\leq N}.
\]
By the similar calculation as \eqref{dispersion:identity} we see the positivity of the symmetric matrix $A_1'$. 
Moreover, we have $\det\tilde{A}(\mbox{\boldmath$0$})=\det A_1'>0$. 
Therefore, we obtain the strict positivity of $\det\tilde{A}(h\xi)$. 

It is easy to see that $\det A(h\xi)$ is a polynomial in $(h|\xi|)^2$ of degree less than or equal to $N+1$ 
and that the coefficient of $(h|\xi|)^{2(N+1)}$ is $\det A_0$. 
By the similar calculation as \eqref{dispersion:identity} we see the positivity of the symmetric matrix $A_0$ 
so that the degree is in fact $N+1$. 
Moreover, $A(\mbox{\boldmath$0$})$ has the zero eigenvalue with an eigenvector $(1,0,\ldots,0)^{\rm T}$, 
so that $\det A(\mbox{\boldmath$0$})=0$. 
Therefore, the term of degree $0$ in $\det A(h\xi)$ vanishes so that $(h|\xi|)^{-2}\det A(h\xi)$ is a 
polynomial in $(h|\xi|)^2$ of degree $N$.
Similarly we see that $\det\tilde{A}(h\xi)$ is a polynomial in $(h|\xi|)^2$ of degree less than or equal to $N$ 
and that the coefficient of $(h|\xi|)^{2N}$ is $\det\tilde{A}_0$, which is also positive, so that the degree 
is in fact $N$. 
\quad$\Box$

\bigskip
Thanks of this proposition and the dispersion relation \eqref{dispersion:dispersion}, the linearized system 
\eqref{dispersion:IK model} is classified into the dispersive system, so that the Isobe--Kakinuma model is 
a nonlinear dispersive system of equations. 
Therefore, we can define the phase speed $c_{IK}(\xi)$ of the plane wave solution to \eqref{dispersion:IK model} 
related to the wave vector $\xi\in\mathbf{R}^n$ by 
\begin{equation}\label{dispersion:phase speed}
c_{IK}(\xi) = \pm\sqrt{gh\frac{(h|\xi|)^{-2}\det A(h\xi)}{\det\tilde{A}(h\xi)}}.
\end{equation}
It follows from Proposition \ref{dispersion:prop 1} that 
\[
\lim_{|\xi|\to+\infty}c_{IK}(\xi)=\pm\sqrt{gh\frac{\det A_0}{\det\tilde{A}_0}},
\]
which is not zero. 
It is not consistent with the linear water waves: 
$\lim_{|\xi|\to+\infty}c_{WW}(\xi)=0$. 
This implies that the Isobe--Kakinuma model \eqref{intro:IK model} cannot give a good approximation 
to the water waves in deep water. 
However, as is shown by the following theorems the Isobe--Kakinuma model gives a very precise approximation 
in the shallow water regime $h|\xi|\ll1$.

\begin{theorem}\label{dispersion:theorem 1}
If we choose $p_i=2i$ $(i=0,1,\ldots,N)$, then $(c_{IK}(\xi))^2$ becomes the $[2N/2N]$ Pad\'e approximant 
of $(c_{WW}(\xi))^2$. 
More precisely, there exists a positive constant $C$ depending only on $N$ such that for any 
$\xi\in\mathbf{R}^n$ and any $h,g>0$ we have 
\[
\biggl|\biggl(\frac{c_{WW}(\xi)}{\sqrt{gh}}\biggr)^2-\biggl(\frac{c_{IK}(\xi)}{\sqrt{gh}}\biggr)^2\biggr|
\leq C(h|\xi|)^{4N+2}.
\]
\end{theorem}

\noindent
{\bf Proof}. \ 
Without loss of generality it is sufficient to consider the case $h=1$ and to show that 
\begin{equation}\label{dispersion:asymptotics}
|\xi|\tanh|\xi|=\frac{\det A(\xi)}{\det\tilde{A}(\xi)}+O(|\xi|^{4N+4})
 \quad\mbox{as}\quad |\xi|\to0.
\end{equation}
Let $\hat{\Phi}(z)=\frac{1}{\cosh|\xi|}\cosh(|\xi|z)$. 
In view of $\hat{\Phi}''(z)=|\xi|^2\hat{\Phi}(z)$, $\hat{\Phi}(1)=1$, $\hat{\Phi}'(1)=|\xi|\tanh|\xi|$, 
and $\hat{\Phi}'(0)=0$, we have 
\begin{equation}\label{dispersion:identity 1}
|\xi|\tanh|\xi|=\int_0^1\bigl\{|\xi|^2(\hat{\Phi}(z))^2+(\hat{\Phi}'(z))^2\bigr\}{\rm d}z.
\end{equation}
By the Taylor series expansion of $\cosh(|\xi|z)$ we have 
\[
\hat{\Phi}(z)=\sum_{i=0}^{2N}\hat{\psi}_iz^{2i}+O(|\xi|^{4N+2}), \quad 
\hat{\Phi}'(z)=\sum_{i=0}^{2N}2i\hat{\psi}_iz^{2i-1}+O(|\xi|^{4N+2})
\]
uniformly with respect to $z\in[0,1]$, where 
\begin{equation}\label{dispersion:psi}
\hat{\psi}_i=\frac{|\xi|^{2i}}{(2i)!\cosh|\xi|}, \quad i=0,1,\ldots.
\end{equation}
Plugging these expansions into \eqref{dispersion:identity 1} we see that 
\[
|\xi|\tanh|\xi|
= \sum_{i,j=0}^N\biggl(\frac{1}{2(i+j)+1}|\xi|^2+\frac{4ij}{2(i+j)-1}\biggr)\hat{\psi}_i\hat{\psi}_j
 + 2R + O(|\xi|^{4N+4}), 
\]
where 
\begin{align*}
R &= \sum_{i=0}^{N-1}\sum_{j=N+1}^{2N}\frac{|\xi|^2}{2(i+j)+1}\hat{\psi}_i\hat{\psi}_j
 + \sum_{i=0}^N\sum_{j=N+1}^{2N}\frac{4ij}{2(i+j)-1}\hat{\psi}_i\hat{\psi}_j \\
&= \sum_{i=0}^{N-1}\sum_{j=N+1}^{2N}\biggl(
 \frac{|\xi|^2}{2(i+j)+1}\hat{\psi}_i\hat{\psi}_j
 +\frac{4(i+1)j}{2(i+j)+1}\hat{\psi}_{i+1}\hat{\psi}_j\biggr) \\
&= \frac{|\xi|}{\cosh^2|\xi|}\sum_{i=0}^{N-1}\sum_{j=N+1}^{2N}
 \frac{|\xi|^{2i+1}}{(2i+1)!}\frac{|\xi|^{2j}}{(2j)!} \\
&= \frac{|\xi|\tanh|\xi|}{\cosh|\xi|}\sum_{j=N+1}^{\infty}\frac{|\xi|^{2j}}{(2j)!}+O(|\xi|^{4N+4}).
\end{align*}
Here we used the explicit form \eqref{dispersion:psi} of $\hat{\psi}_i$. 
Therefore, we obtain 
\begin{align}\label{dispersion:identity 2}
& (|\xi|\tanh|\xi|)\biggl(1-\frac{2}{\cosh|\xi|}\sum_{j=N+1}^{\infty}\frac{|\xi|^{2j}}{(2j)!}\biggr) \\
&= \sum_{i,j=0}^N\biggl(\frac{1}{2(i+j)+1}|\xi|^2+\frac{4ij}{2(i+j)-1}\biggr)\hat{\psi}_i\hat{\psi}_j
 + O(|\xi|^{4N+4}). \nonumber
\end{align}

Up to now we did not use the specific choice of $p_i$. 
From now on, we use the advantage of the choice $p_i=2i$. 
Then, we can rewrite the above identity \eqref{dispersion:identity 2} as 
\begin{equation}\label{dispersion:identity 3}
 (|\xi|\tanh|\xi|)\biggl(1-\frac{2}{\cosh|\xi|}\sum_{j=N+1}^{\infty}\frac{|\xi|^{2j}}{(2j)!}\biggr) \\
= (A(\xi)\hat{\mbox{\boldmath$\psi$}})\cdot\hat{\mbox{\boldmath$\psi$}} + O(|\xi|^{4N+4}),
\end{equation}
where $\hat{\mbox{\boldmath$\psi$}}=(\hat{\psi}_0,\ldots,\hat{\psi}_N)^{\rm T}$. 
On the other hand, we consider the following equations for $\hat{\eta}$ and 
$\hat{\mbox{\boldmath$\phi$}}=(\hat{\phi}_0,\ldots,\hat{\phi}_N)^{\rm T}$. 
\begin{equation}\label{dispersion:equation}
\left(
 \begin{array}{cc}
  0 & \mbox{\boldmath$1$}^{\rm T} \\
  -\mbox{\boldmath$1$} & A(\xi) 
 \end{array}
\right)
\left(
 \begin{array}{c}
  \hat{\eta} \\
  \hat{\mbox{\boldmath$\phi$}}
 \end{array}
\right)
=
\left(
 \begin{array}{c}
  1 \\
  \mbox{\boldmath$0$}
 \end{array}
\right).
\end{equation}
Since the coefficient matrix $\tilde{A}(\xi)$ is nonsingular, we have a unique solution 
$(\hat{\eta},\hat{\mbox{\boldmath$\phi$}})$ of these equations. 
Moreover, we see that 
\begin{equation}\label{dispersion:identity 4}
(A(\xi)\hat{\mbox{\boldmath$\phi$}})\cdot\hat{\mbox{\boldmath$\phi$}}
=\hat{\eta}=\frac{\det A(\xi)}{\det\tilde{A}(\xi)}. 
\end{equation}
In view of this and \eqref{dispersion:identity 3}, it is sufficient to compare 
$\hat{\mbox{\boldmath$\psi$}}$ with $\hat{\mbox{\boldmath$\phi$}}$. 
Let $A(\xi)=(\mbox{\boldmath$a$}_0,\ldots,\mbox{\boldmath$a$}_N)$. 
Then, by using the explicit form \eqref{dispersion:psi} of $\hat{\psi}_i$ we see that 
\[
\mbox{\boldmath$a$}_i\cdot\hat{\mbox{\boldmath$\psi$}}
= \frac{1}{\cosh|\xi|}\biggl(\sum_{j=0}^{N-1}\frac{|\xi|^{2j+2}}{(2j+1)!}
 +\frac{1}{2(i+N)+1}\frac{|\xi|^{2N+2}}{(2N)!}\biggr)
\]
for $i=0,1,\ldots,N$, and that 
\[
\mbox{\boldmath$1$}\cdot\hat{\mbox{\boldmath$\psi$}}
=\frac{1}{\cosh|\xi|}\sum_{i=0}^N\frac{|\xi|^{2i}}{(2i)!}
=1+O(|\xi|^{2N+2}).
\]
Therefore, if we put $\hat{\zeta}=\frac{1}{\cosh|\xi|}\sum_{j=0}^{N-1}\frac{|\xi|^{2j+2}}{(2j+1)!}$, 
then we have 
\begin{equation}\label{dispersion:relation}
\left(
 \begin{array}{cc}
  0 & \mbox{\boldmath$1$}^{\rm T} \\
  -\mbox{\boldmath$1$} & A(\xi) 
 \end{array}
\right)
\left(
 \begin{array}{c}
  \hat{\zeta} \\
  \hat{\mbox{\boldmath$\psi$}}
 \end{array}
\right)
=
\left(
 \begin{array}{c}
  1 \\
  \mbox{\boldmath$0$}
 \end{array}
\right)+O(|\xi|^{2N+2}).
\end{equation}
By Proposition \ref{dispersion:prop 1} the determinant of the coefficient matrix is strictly positive, 
so that this together with \eqref{dispersion:equation} implies 
$\hat{\mbox{\boldmath$\psi$}}=\hat{\mbox{\boldmath$\phi$}}+O(|\xi|^{2N+2})$. 
Moreover, it holds that 
\[
(A(\xi)\hat{\mbox{\boldmath$\phi$}})\cdot\hat{\mbox{\boldmath$\psi$}}
= \frac{\det A(\xi)}{\det\tilde{A}(\xi)}\frac{1}{\cosh|\xi|}\sum_{i=0}^N\frac{|\xi|^{2i}}{(2i)!},
\]
which together with \eqref{dispersion:identity 4} implies 
\[
(A(\xi)\hat{\mbox{\boldmath$\phi$}})\cdot(\hat{\mbox{\boldmath$\phi$}}-\hat{\mbox{\boldmath$\psi$}})
= \frac{\det A(\xi)}{\det\tilde{A}(\xi)}\frac{1}{\cosh|\xi|}\sum_{i=N+1}^{\infty}\frac{|\xi|^{2i}}{(2i)!},
\]
so that 
\begin{align}\label{dispersion:relation 2}
(A(\xi)\hat{\mbox{\boldmath$\psi$}})\cdot\hat{\mbox{\boldmath$\psi$}}
&= (A(\xi)\hat{\mbox{\boldmath$\phi$}})\cdot\hat{\mbox{\boldmath$\phi$}}
 -2(A(\xi)\hat{\mbox{\boldmath$\phi$}})\cdot(\hat{\mbox{\boldmath$\phi$}}-\hat{\mbox{\boldmath$\psi$}})
 +(A(\xi)(\hat{\mbox{\boldmath$\phi$}}-\hat{\mbox{\boldmath$\psi$}}))
  \cdot(\hat{\mbox{\boldmath$\phi$}}-\hat{\mbox{\boldmath$\psi$}}) \\
&= \frac{\det A(\xi)}{\det\tilde{A}(\xi)}
 \biggl(1-\frac{2}{\cosh|\xi|}\sum_{i=N+1}^{\infty}\frac{|\xi|^{2i}}{(2i)!}\biggr)
 +O(|\xi|^{4N+4}). \nonumber 
\end{align}
Plugging this into \eqref{dispersion:identity 3} we obtain the desired relation \eqref{dispersion:asymptotics}. 
\quad$\Box$

\bigskip
As was shown by J. Boussinesq \cite{Boussinesq1872}, in the case of the flat bottom the velocity 
potential $\Phi$ for the water wave problem can be expanded in a Taylor series with respect to 
the vertical spatial variable $z$ around the bottom $z=-h$ as 
$$
\Phi(x,z,t)=\sum_{i=0}^{\infty}\frac{(z+h)^{2i}}{(2i)!}(-\Delta)^i\phi_0(x,t),
$$
where $\phi_0$ is the trace of the velocity potential $\Phi$ on the bottom. 
Therefore, it is natural to choose $p_i=2i$ in such a case. 
However, in the case of variable bottom topographies the Taylor series with respect to $z$ around 
the bottom $z=-h+b(x)$ contains terms of odd degree too, so that the choice $p_i=i$ would be important 
to such cases. 
If we choose $p_i=i$, we do not have such a beautiful result as Theorem \ref{dispersion:theorem 1} but 
we still have the following theorem, 
which asserts that the terms of odd degree in the approximation of the velocity potential $\Phi$ 
do not affect the precision of the linear dispersion relation in the shallow water regime $h|\xi|\ll1$.

\begin{theorem}\label{dispersion:theorem 2}
If we choose $p_i=i$ $(i=0,1,\ldots,N)$, then for any $\xi\in\mathbf{R}^n$ and any $h,g>0$ we have 
\[
\biggl|\biggl(\frac{c_{WW}(\xi)}{\sqrt{gh}}\biggr)^2-\biggl(\frac{c_{IK}(\xi)}{\sqrt{gh}}\biggr)^2\biggr|
\leq C(h|\xi|)^{4[N/2]+2},
\]
where $C$ is a positive constant depending only on $N$ and $[N/2]$ is the integer part of $N/2$. 
\end{theorem}

\noindent
{\bf Proof}. 
We will give the proof in the case of even integer $N$. 
As in the proof of the previous theorem, it is sufficient to show that 
\begin{equation}\label{dispersion:asymptotics 2}
|\xi|\tanh|\xi|=\frac{\det A(\xi)}{\det\tilde{A}(\xi)}+O(|\xi|^{4[N/2]+4})
 \quad\mbox{as}\quad |\xi|\to0.
\end{equation}
In place of \eqref{dispersion:psi}, we define $\hat{\psi}_i$ by 
\[
\hat{\psi}_i=
\left\{
 \begin{array}{ll}
  \displaystyle
  \frac{|\xi|^i}{i!\cosh|\xi|} & \mbox{if $i$ is even}, \\[0.5ex]
  0 & \mbox{if $i$ is odd}
 \end{array}
\right.
\]
and put $\hat{\mbox{\boldmath$\psi$}}=(\hat{\psi}_0,\ldots,\hat{\psi}_N)^{\rm T}$. 
Then, we have 
\[
(A(\xi)\hat{\mbox{\boldmath$\psi$}})\cdot\hat{\mbox{\boldmath$\psi$}}
=\sum_{i,j=0}^{[N/2]}\biggl(\frac{1}{2(i+j)+1}|\xi|^2+\frac{4ij}{2(i+j)-1}\biggr)\hat{\psi}_{2i}\hat{\psi}_{2j},
\]
so that we can rewrite the identity \eqref{dispersion:identity 2} as 
\begin{equation}\label{dispersion:identity 5}
 (|\xi|\tanh|\xi|)\biggl(1-\frac{2}{\cosh|\xi|}\sum_{j=[N/2]+1}^{\infty}\frac{|\xi|^{2j}}{(2j)!}\biggr) \\
= (A(\xi)\hat{\mbox{\boldmath$\psi$}})\cdot\hat{\mbox{\boldmath$\psi$}} + O(|\xi|^{4[N/2]+4}).
\end{equation}
Let $(\hat{\eta},\hat{\mbox{\boldmath$\phi$}})$ be the solution of \eqref{dispersion:equation} as before. 
Then, we have the identity \eqref{dispersion:identity 4}. 
If we put $\hat{\zeta}=\frac{1}{\cosh|\xi|}\sum_{j=0}^{[N/2]-1}\frac{|\xi|^{2j+2}}{(2j+1)!}$, 
then in place of \eqref{dispersion:relation} we have 
\[
\left(
 \begin{array}{cc}
  0 & \mbox{\boldmath$1$}^{\rm T} \\
  -\mbox{\boldmath$1$} & A(\xi) 
 \end{array}
\right)
\left(
 \begin{array}{c}
  \hat{\zeta} \\
  \hat{\mbox{\boldmath$\psi$}}
 \end{array}
\right)
=
\left(
 \begin{array}{c}
  1 \\
  \mbox{\boldmath$0$}
 \end{array}
\right)+O(|\xi|^{2[N/2]+2}).
\]
This together with \eqref{dispersion:equation} implies 
$\hat{\mbox{\boldmath$\psi$}}=\hat{\mbox{\boldmath$\phi$}}+O(|\xi|^{2[N/2]+2})$. 
Therefore, in place of \eqref{dispersion:relation 2} we have 
\[
(A(\xi)\hat{\mbox{\boldmath$\psi$}})\cdot\hat{\mbox{\boldmath$\psi$}}
= \frac{\det A(\xi)}{\det\tilde{A}(\xi)}
 \biggl(1-\frac{2}{\cosh|\xi|}\sum_{i=[N/2]+1}^{\infty}\frac{|\xi|^{2i}}{(2i)!}\biggr)
 +O(|\xi|^{4[N/2]+4}). 
\]
Plugging this into \eqref{dispersion:identity 5} we obtain the desired relation \eqref{dispersion:asymptotics 2}. 
The case of odd integer $N$ can be proved in the same way, so we omit it. 
\quad$\Box$

%------------------------------------------------------------------------------
\section{Analysis of a linearized system}
\label{section:linear}
\setcounter{equation}{0}
\setcounter{theorem}{0}
In this section we consider the linearized equations of the Isobe--Kakinuma model \eqref{intro:IK model} 
around an arbitrary flow $(\eta,\phi_0,\ldots,\phi_N)$, which is assumed to be sufficiently smooth. 
The hypersurface $t=0$ is still characteristic for the linearized equations. 
We will transform the equations to a symmetric positive system of partial differential equations for which 
the hypersurface $t=0$ is noncharacteristic and give an idea to derive a priori estimates for the solution 
to the nonlinear equations. 

We introduce second order differential operators $L_{ij}=L_{ij}(H,b)$ $(i,j=0,1,\ldots,N)$ depending on 
the water depth $H$ and the bottom topography $b$ by 
\begin{align}\label{linear:L}
L_{ij}\psi_j
&= -\nabla\cdot\biggl(
   \frac{1}{p_i+p_j+1}H^{p_i+p_j+1}\nabla\psi_j
   -\frac{p_j}{p_i+p_j}H^{p_i+p_j}\psi_j\nabla b\biggr) \\[0.5ex]
&\quad\,
  -\frac{p_i}{p_i+p_j}H^{p_i+p_j}\nabla b\cdot\nabla\psi_j
   +\frac{p_ip_j}{p_i+p_j-1}H^{p_i+p_j-1}(1+|\nabla b|^2)\psi_j. \nonumber
\end{align}
Then, we have $L_{ij}^*=L_{ji}$, where $L_{ij}^*$ is the adjoint operator of $L_{ij}$ in $L^2(\mathbf{R}^n)$. 
In addition to the function $a$ defined by \eqref{intro:a} we introduce the functions $\mbox{\boldmath$u$}$ 
and $w$ by 
\begin{equation}\label{linear:uw}
\mbox{\boldmath$u$}=\sum_{i=0}^N(H^{p_i}\nabla\phi_i-p_iH^{p_i-1}\phi_i\nabla b), \quad
w=\sum_{i=0}^Np_iH^{p_i-1}\phi_i.
\end{equation}
Since $\mbox{\boldmath$u$}=\nabla\Phi^{\mbox{\rm\tiny app}}|_{z=\eta}$ and 
$w=\partial_z\Phi^{\mbox{\rm\tiny app}}|_{z=\eta}$, where $\Phi^{\mbox{\rm\tiny app}}$ is the approximate 
velocity potential defined by \eqref{intro:approximation}, $\mbox{\boldmath$u$}$ and $w$ represent 
approximately the horizontal and the vertical components of the velocity field on the water surface, 
respectively. 
Then, the Isobe--Kakinuma model \eqref{intro:IK model} can be written simply as 
\[
\left\{
 \begin{array}{l}
  \displaystyle
  H^{p_i}\partial_t\eta-\sum_{j=0}^NL_{ij}\phi_j=0 \qquad\mbox{for}\quad i=0,1,\ldots,N, \\
  \displaystyle
  \sum_{j=0}^NH^{p_j}\partial_t\phi_j+g\eta+\frac12(|\mbox{\boldmath$u$}|^2+w^2)=0.
 \end{array}
\right.
\]
Now, let us linearize the above equations around $(\eta,\phi_0,\ldots,\phi_N)$. 
We denote by $(\zeta,\psi_0,\ldots,\psi_N)$ the variation from $(\eta,\phi_0,\ldots,\phi_N)$. 
After a tedious but straightforward calculation we obtain the linearized equations 
\begin{equation}\label{linear:linearized equations}
\left\{
 \begin{array}{l}
  \displaystyle
  H^{p_i}(\partial_t\zeta+\nabla\cdot(\mbox{\boldmath$u$}\zeta))
   +p_iH^{p_i-1}(\partial_t\eta+\mbox{\boldmath$u$}\cdot\nabla\eta-w)\zeta
   -\sum_{j=0}^NL_{ij}\psi_j=f_i \\
  \makebox[20em]{}\mbox{for}\quad i=0,1,\ldots,N, \\
  \displaystyle
  \sum_{j=0}^N\Bigl\{H^{p_j}(\partial_t\psi_j+\mbox{\boldmath$u$}\cdot\nabla\psi_j)
   -p_jH^{p_j-1}(\mbox{\boldmath$u$}\cdot\nabla b-w)\psi_j\Bigr\}+a\zeta=f_{N+1},
 \end{array}
\right.
\end{equation}
where $f_0,\ldots,f_{N+1}$ are given functions.

The hypersurface $t=0$ is still characteristic for these linearized equations. 
In fact, by eliminating the time derivative $\partial_t\zeta$ from these equations we have 
\[
\sum_{j=0}^N(L_{ij}-H^{p_i}L_{0j})\psi_j
=p_iH^{p_i-1}(\partial_t\eta+\mbox{\boldmath$u$}\cdot\nabla\eta-w)\zeta+H^{p_i}f_0-f_i
\]
for $i=1,\ldots,N$. 
Now, we differentiate this with respect to the time $t$ and use the first equation in 
\eqref{linear:linearized equations} to eliminate the time derivative $\partial_t\zeta$. 
Then, we obtain 
\begin{align}\label{linear:reduced 1}
\sum_{j=0}^N(L_{ij}-H^{p_i}L_{0j})\partial_t\psi_j
&= p_iH^{p_i-1}(\partial_t\eta+\mbox{\boldmath$u$}\cdot\nabla\eta-w)
 \biggl(\sum_{j=0}^NL_{0j}\psi_j-\nabla\cdot(\mbox{\boldmath$u$}\zeta)\biggr) \\
&\quad\;
 -\sum_{j=0}^N[\partial_t,L_{ij}-H^{p_i}L_{0j}]\psi_j
 +[\partial_t,p_iH^{p_i-1}(\partial_t\eta+\mbox{\boldmath$u$}\cdot\nabla\eta-w)]\zeta \nonumber \\
&\quad\;
 +\partial_t(H^{p_i}f_0-f_i)+p_iH^{p_i-1}(\partial_t\eta+\mbox{\boldmath$u$}\cdot\nabla\eta-w)f_0.
 \nonumber
\end{align}
On the other hand, it follows from the second equation in \eqref{linear:linearized equations} that 
\begin{equation}\label{linear:reduced 2}
\sum_{j=0}^NH^{p_j}\partial_t\psi_j
=-\sum_{j=0}^N\Bigl\{H^{p_j}\mbox{\boldmath$u$}\cdot\nabla\psi_j
   -p_jH^{p_j-1}(\mbox{\boldmath$u$}\cdot\nabla b-w)\psi_j\Bigr\}-a\zeta+f_{N+1}.
\end{equation}
Here we note that the right hand sides of \eqref{linear:reduced 1} and \eqref{linear:reduced 2} do not include any 
time derivatives of $(\zeta,\psi_0,\ldots,\psi_N)$. 
In view of these equations we introduce linear operators $\mathscr{L}_i=\mathscr{L}_i(H,b)$ $(i=0,1,\ldots,N)$ 
depending on the water depth $H$ and the bottom topography $b$ and acting on 
$\mbox{\boldmath$\varphi$}=(\varphi_0,\ldots,\varphi_N)^{\rm T}$ by 
\begin{equation}\label{linear:L2}
\mathscr{L}_0\mbox{\boldmath$\varphi$}=\sum_{j=0}^NH^{p_j}\varphi_j, \qquad
\mathscr{L}_i\mbox{\boldmath$\varphi$}=\sum_{j=0}^N(L_{ij}-H^{p_i}L_{0j})\varphi_j
 \quad\mbox{for}\quad i=1,\ldots,N, 
\end{equation}
and put $\mathscr{L}\mbox{\boldmath$\varphi$}=(\mathscr{L}_0\mbox{\boldmath$\varphi$},\ldots,
\mathscr{L}_N\mbox{\boldmath$\varphi$})^{\rm T}$. 
For a given $\mbox{\boldmath$F$}=(F_0,\ldots,F_N)^{\rm T}$ we consider the equation 
\begin{equation}\label{linear:elliptic 1}
\mathscr{L}\mbox{\boldmath$\varphi$}=\mbox{\boldmath$F$}.
\end{equation}
Once we show the solvability of $\mbox{\boldmath$\varphi$}$ of this equation, we could express the 
time derivatives $\partial_t\psi_i$ $(i=0,1,\ldots,N)$ in terms of the spatial derivatives so that 
we could avoid the difficulty arising from the fact that the hypersurface $t=0$ is characteristic for 
the linearized equations \eqref{linear:linearized equations}.

Let $\mbox{\boldmath$\varphi$}$ be a solution of \eqref{linear:elliptic 1}. 
It follows from the first component of \eqref{linear:elliptic 1} that 
\begin{equation}\label{linear:reduction}
\varphi_0=F_0-\sum_{j=1}^NH^{p_j}\varphi_j.
\end{equation}
Plugging this into the other components of \eqref{linear:elliptic 1} we obtain 
\begin{equation}\label{linear:elliptic 2}
P_i\mbox{\boldmath$\varphi$}'= F_i-(L_{i0}-H^{p_i}L_{00})F_0 \quad\mbox{for}\quad i=1,\ldots,N,
\end{equation}
where $\mbox{\boldmath$\varphi$}'=(\varphi_1,\ldots,\varphi_N)^{\rm T}$ and 
$P_j=P_j(H,b)$ $(j=1,\ldots,N)$ are second order differential operators defined by 

\begin{equation}\label{linear:P}
P_i\mbox{\boldmath$\varphi$}'=
\sum_{j=1}^N\Bigl\{(L_{ij}-H^{p_i}L_{0j})\varphi_j-(L_{i0}-H^{p_i}L_{00})(H^{p_j}\varphi_j)\Bigr\}. 
\end{equation}
We further introduce the operator $P\mbox{\boldmath$\varphi$}'=(P_1\mbox{\boldmath$\varphi$}',\ldots,
P_N\mbox{\boldmath$\varphi$}')^{\rm T}$. 
Since $L_{ij}^*=L_{ji}$, we see easily that $P$ is symmetric in $L^2(\mathbf{R}^n)$. 
Moreover, $P$ is positive in $L^2(\mathbf{R}^n)$ as shown in the following lemma.

\begin{lemma}\label{linear:lemma 1}
Let $c_0, c_1$ be positive constants. 
There exists a positive constant $C=C(c_0,c_1)>0$ depending only on $c_0$ and $c_1$ such that 
if $H,\nabla b\in L^{\infty}(\mathbf{R}^n)$ satisfy $H(x)\geq c_0$ and $|\nabla b(x)|\leq c_1$, 
then we have 
\[
(P\mbox{\boldmath$\varphi$}',\mbox{\boldmath$\varphi$}')_{L^2}\geq C^{-1}\|\mbox{\boldmath$\varphi$}'\|_1^2. 
\]
\end{lemma}

\noindent
{\bf Proof}. \ 
Introducing $\varphi_0=-\sum_{j=1}^NH^{p_j}\varphi_j$, we have 
\begin{equation}\label{linear:equality 1}
(P\mbox{\boldmath$\varphi$}',\mbox{\boldmath$\varphi$}')_{L^2}
=\sum_{i,j=0}^N(L_{ij}\varphi_j,\varphi_i)_{L^2}.
\end{equation}
Although this equality can be derived by direct calculation, it is also derived by the following argument 
which help us to understand the structure of the equations. 
Put $P\mbox{\boldmath$\varphi$}'=(G_1,\ldots,G_N)^{\rm T}$. 
Then, we have 
\[
\sum_{j=0}^NH^{p_j}\varphi_j=0, \qquad
\sum_{j=0}^N(L_{ij}-H^{p_i}L_{0j})\varphi_j=G_i \quad\mbox{for}\quad i=1,\ldots,N.
\]
If we further introduce $\zeta=\sum_{j=0}^NL_{0j}\varphi_j$, then the above equations can be written as 
\[
\left(
 \begin{array}{cc}
  0 & \mbox{\boldmath$l$}^{\rm T} \\
  -\mbox{\boldmath$l$} & L
 \end{array}
\right)
\left(
 \begin{array}{c}
  \zeta \\
  \mbox{\boldmath$\varphi$}
 \end{array}
\right)
=
\left(
 \begin{array}{c}
  0 \\
  \mbox{\boldmath$G$}
 \end{array}
\right),
\]
where $\mbox{\boldmath$\varphi$}=(\varphi_0,\ldots,\varphi_N)^{\rm T}$, 
$\mbox{\boldmath$G$}=(0,G_1,\ldots,G_N)^{\rm T}$, $L=(L_{ij})_{0\leq i,j\leq N}$, and 
\begin{equation}\label{linear:l}
\mbox{\boldmath$l$}=\mbox{\boldmath$l$}(H)=(H^{p_0},\ldots,H^{p_N})^{\rm T}.
\end{equation}
By taking the $L^2$-inner product of the above equation with $(\zeta,\mbox{\boldmath$\varphi$}^{\rm T})^{\rm T}$, 
we obtain $(L\mbox{\boldmath$\varphi$},\mbox{\boldmath$\varphi$})_{L^2}
=(\mbox{\boldmath$G$},\mbox{\boldmath$\varphi$})_{L^2}
=(P\mbox{\boldmath$\varphi$}',\mbox{\boldmath$\varphi$}')_{L^2}$, which implies \eqref{linear:equality 1}. 

By direct calculation we have 
\begin{align}\label{linear:equivalence}
I &:= \int_{\mathbf{R}^n}\!{\rm d}x\!\int_0^{H}\biggl\{
  \biggl|\sum_{i=0}^N(z^{p_i}\nabla\varphi_i-p_iz^{p_i-1}\varphi_i\nabla b)\biggr|^2
 +\biggl(\sum_{i=0}^Np_iz^{p_i-1}\varphi_i\biggr)^2\biggr\}{\rm d}z \\
&= \sum_{i,j=0}^N(L_{ij}\varphi_i,\varphi_j)_{L^2}. \nonumber 
\end{align}
In the case $p_i=2i$ $(i=0,1,\ldots,N)$, $\{z^{p_i},z^{p_i-1}\}_{0\leq i\leq N}$ are linearly independent 
so that 
\begin{align*}
I &\simeq \int_{\mathbf{R}^n}\!{\rm d}x\!\int_0^{H}\sum_{i=0}^N\Bigl\{
 \bigl(z^{2p_i}|\nabla\varphi_i|^2+p_i^2z^{2p_i-2}|\nabla b|^2\varphi_i^2\bigr)
 +p_i^2z^{2p_i-2}\varphi_i^2\Bigr\}{\rm d}z \\
&\simeq \int_{\mathbf{R}^n}\sum_{i=0}^N\Bigl\{
 H^{2p_i+1}|\nabla\varphi_i|^2+p_i^2H^{2p_i-1}(1+|\nabla b|^2)\varphi_i^2\Bigr\}{\rm d}x,
\end{align*}
which together with \eqref{linear:equality 1} and \eqref{linear:equivalence} implies the desired estimate 
in that case. 
We remark that the constant $C$ can be taken independent of $c_1$. 
In the case $p_i=i$ $(i=0,1,\ldots,N)$, we see that 
\begin{align*}
I &= \int_{\mathbf{R}^n}\!{\rm d}x\!\int_0^{H}\biggl\{
  \biggl|\sum_{i=0}^{N-1}z^{i}\bigl(\nabla\varphi_i-(i+1)\varphi_{i+1}\nabla b\bigr)
   +z^N\nabla\varphi_N\biggr|^2
 +\biggl(\sum_{i=1}^Niz^{i-1}\varphi_i\biggr)^2\biggr\}{\rm d}z \\
&\simeq \int_{\mathbf{R}^n}\biggl\{\sum_{i=0}^{N-1}H^{2i+1}|\nabla\varphi_i-(i+1)\varphi_{i+1}\nabla b|^2
 +H^{N+1}|\nabla\varphi_N|^2+\sum_{i=1}^NH^{2i-1}\varphi_i^2\biggr\}{\rm d}z \\
&\gtrsim \sum_{i=0}^{N-1}\|\nabla\varphi_i-(i+1)\varphi_{i+1}\nabla b\|^2
 +\|\nabla\varphi_N\|^2+\sum_{i=1}^N\|\varphi_i\|^2
\end{align*}
which together with \eqref{linear:equality 1} and \eqref{linear:equivalence} implies the desired estimate 
in that case. 
The other cases can be treated in the same way, so we omit it. 
\quad$\Box$

\bigskip
By this lemma, the explicit expression \eqref{linear:P} of the operator $P$ (see also \eqref{linear:L}), 
and the standard theory of elliptic partial differential equations, we can obtain the following lemma.

\begin{lemma}\label{linear:lemma 2}
Let $h,c_0,M$ be positive constants and $m$ an integer such that $m>n/2+1$. 
There exists a positive constant $C=C(h,c_0,M)$ such that if $\eta$ and $b$ satisfy 
\[
\left\{
 \begin{array}{l}
  \|\eta\|_m+\|b\|_{W^{m,\infty}} \leq M, \\[0.5ex]
  c_0\leq H(x)=h+\eta(x)-b(x) \quad\mbox{for}\quad x\in\mathbf{R}^n,
 \end{array}
\right.
\]
then for $1\leq k\leq m$ we have 
\[
\|P^{-1}\mbox{\boldmath$G$}'\|_k \leq C\|\mbox{\boldmath$G$}'\|_{k-2}.
\]
\end{lemma}

We proceed to consider equation \eqref{linear:elliptic 1}. 
Thanks of this lemma, we see that for a given $\mbox{\boldmath$F$}$ there exists a unique solution 
$\mbox{\boldmath$\varphi$}'=(\varphi_1,\ldots,\varphi_N)^{\rm T}$ of \eqref{linear:elliptic 2}. 
If we define $\varphi_0$ by \eqref{linear:reduction}, then 
$\mbox{\boldmath$\varphi$}=(\varphi_0,\ldots,\varphi_N)^{\rm T}$ is a solution of \eqref{linear:elliptic 1}. 
More precisely, we have the following lemma.

\begin{lemma}\label{linear:lemma 3}
Under the hypothesis of Lemma {\rm \ref{linear:lemma 2}}, for any 
$\mbox{\boldmath$F$}=(F_0,\ldots,F_N)^{\rm T}$ satisfying $\nabla F_0\in H^{k-1}$ and 
$(F_1,\ldots,F_N)\in H^{k-2}$ with $1\leq k\leq m$ there exists a unique solution 
$\mbox{\boldmath$\varphi$}=(\varphi_0,\ldots,\varphi_N)^{\rm T}$ of \eqref{linear:elliptic 1} satisfying 
\[
\|\nabla\varphi_0\|_{k-1}+\|(\varphi_1,\ldots,\varphi_N)\|_k
\leq C(\|\nabla F_0\|_{k-1}+\|(F_1,\ldots,F_N)\|_{k-2}),
\]
where $C=C(h,c_0,M)>0$. 
If, in addition, $F_0\in L^2(\mathbf{R}^n)$, then we have 
\[
\|\mbox{\boldmath$\varphi$}\|_k \leq C(\|F_0\|_k+\|(F_1,\ldots,F_N)\|_{k-2}).
\]
\end{lemma}

\medskip
Now, we have established the solvability of $\mbox{\boldmath$\varphi$}$ to equation \eqref{linear:elliptic 1}. 
We go back to consider the linearized equations in \eqref{linear:linearized equations}. 
In the following of this section we denote lower order terms and terms related with the given functions 
$f_0,\ldots,f_{N+1}$ by the same symbol $LOT$, which may change from line to line. 
We introduce a symmetric matrix $A(H)$ depending on the water depth $H$ by 
\begin{equation}\label{linear:A}
A(H)=(a_{ij}(H))_{0\leq i,j\leq N}, \quad a_{ij}(H)=\frac{1}{p_i+p_j+1}H^{p_i+p_j+1}.
\end{equation}
We also use a matrix $\tilde{A}(H)$ defined by 
\begin{equation}\label{linear:tildeA}
\tilde{A}(H)=
\left(
 \begin{array}{cc}
  0 & \mbox{\boldmath$l$}(H)^{\rm T} \\
  -\mbox{\boldmath$l$}(H) & A(H)
 \end{array}
\right),
\end{equation}
where $\mbox{\boldmath$l$}(H)$ is defined by \eqref{linear:l}. 
Since $L_{ij}(H,b)=-a_{ij}(H)\Delta+LOT$, it follows from \eqref{linear:reduced 1} and \eqref{linear:reduced 2} 
that 
\[
\left\{
 \begin{array}{l}
  \displaystyle
  \sum_{j=0}^N(a_{ij}(H)-H^{p_i}a_{0j}(H))\Delta\partial_t\psi_j=\nabla\cdot(LOT)+LOT
   \quad\mbox{for}\quad i=1,\ldots,N, \\
  \displaystyle
  \sum_{j=0}^NH^{p_j}\Delta\partial_t\psi_j
   +\sum_{j=0}^NH^{p_j}(\mbox{\boldmath$u$}\cdot\nabla)\Delta\psi_j+\Delta(a\zeta)=\nabla\cdot(LOT)+LOT.
 \end{array}
\right.
\]
We can rewrite these equations in a matrix form as 
\begin{align}\label{linear:transformation 1}
\tilde{A}(H)
\left(
 \begin{array}{c}
  \sum_{j=0}^Na_{0j}(H)\Delta\partial_t\psi_j \\
  \Delta\partial_t\mbox{\boldmath$\psi$}
 \end{array}
\right)
+ &
\left(
 \begin{array}{c}
  \sum_{j=0}^NH^{p_j}(\mbox{\boldmath$u$}\cdot\nabla)\Delta\psi_j+\Delta(a\zeta) \\
  \mbox{\boldmath$0$}
 \end{array}
\right) \\
&=
  \nabla\cdot(LOT)+LOT, \nonumber
\end{align}
where $\mbox{\boldmath$\psi$}=(\psi_0,\ldots,\psi_N)^{\rm T}$. 
Since $\det\tilde{A}(H)=H^{2\sum_{i=0}^Np_i+N}\det\tilde{A}_0$, where $\tilde{A}_0$ is the matrix defined 
in Section \ref{section:dispersion}, the matrix $\tilde{A}(H)$ is nonsingular and its inverse matrix can 
be written as 
\begin{equation}\label{linear:inverse}
\tilde{A}(H)^{-1}=
\left(
 \begin{array}{cc}
  q(H) & \mbox{\boldmath$q$}(H)^{\rm T} \\
  -\mbox{\boldmath$q$}(H) & Q(H)
 \end{array}
\right)
\end{equation}
with a symmetric matrix $Q(H)$. 
If it causes no confusion, we omit the dependence of $H$ in the notation. 
Then, it holds that 
\begin{equation}\label{linear:q}
\mbox{\boldmath$l$}\cdot\mbox{\boldmath$q$}=-1, \quad
A\mbox{\boldmath$q$}=-q\mbox{\boldmath$l$}, \quad
q=\frac{\det A}{\det\tilde{A}}=H\frac{\det A_0}{\det\tilde{A}_0}.
\end{equation}
Therefore, it follows from \eqref{linear:transformation 1} that 
\begin{equation}\label{linear:transformation 2}
-A\Delta\partial_t\mbox{\boldmath$\psi$}
= q\mbox{\boldmath$l$}\bigl\{\mbox{\boldmath$l$}^{\rm T}(\mbox{\boldmath$u$}\cdot\nabla)\Delta\mbox{\boldmath$\psi$}
 +\Delta(a\zeta)\bigr\}+\nabla\cdot(LOT)+LOT. 
\end{equation}
On the other hand, it follows from the first equation in \eqref{linear:linearized equations} that 
\[
\mbox{\boldmath$l$}(\partial_t\zeta+\mbox{\boldmath$u$}\cdot\nabla\zeta)+A\Delta\mbox{\boldmath$\psi$}=LOT.
\]
Taking the Euclidean inner product of this equation with $-a\mbox{\boldmath$q$}$ and using the relations in 
\eqref{linear:q} we obtain 
\[
a(\partial_t\zeta+\mbox{\boldmath$u$}\cdot\nabla\zeta)
 +a\Delta(q\mbox{\boldmath$l$}^{\rm T}\mbox{\boldmath$\psi$})=LOT.
\]
We can write this equation and \eqref{linear:transformation 2} in a matrix form as 
\begin{align}\label{linear:positive system}
\left(
 \begin{array}{cc}
  a & \mbox{\boldmath$0$}^{\rm T} \\
  \mbox{\boldmath$0$} & -A\Delta
 \end{array}
\right)
\partial_t
\left(
 \begin{array}{c}
  \zeta \\
  \mbox{\boldmath$\psi$}
 \end{array}
\right)
+ &
\left(
 \begin{array}{cc}
  a\mbox{\boldmath$u$}\cdot\nabla & a\Delta(q\mbox{\boldmath$l$}^{\rm T}\,\cdot\,) \\
  -q\mbox{\boldmath$l$}\Delta(a\,\cdot\,) & 
   -q\mbox{\boldmath$l$}\mbox{\boldmath$l$}^{\rm T}(\mbox{\boldmath$u$}\cdot\nabla)\Delta
 \end{array}
\right)
\left(
 \begin{array}{c}
  \zeta \\
  \mbox{\boldmath$\psi$}
 \end{array}
\right) \\
&=
\left(
 \begin{array}{c}
  LOT \\
  \nabla\cdot(LOT)+LOT
 \end{array}
\right). \nonumber
\end{align}
Here the operator in the second term in the left-hand side is skew-symmetric in $L^2(\mathbf{R}^n)$ 
modulo lower order terms. 
Since the matrix $A=A(H)$ can be written as 
\[
A(H)=H{\rm diag}(H^{p_0},\ldots,H^{p_N})A_0{\rm diag}(H^{p_0},\ldots,H^{p_N})
\]
and we have shown the positivity of the matrix $A_0$ in Proposition \ref{dispersion:prop 1}, 
$A(H)$ is strictly positive under the strict positivity of the water depth $H$. 
Therefore, under the sign condition $a>0$ \eqref{linear:positive system} forms a symmetric positive system, 
so that the corresponding energy function is defined by 
\[
E(\zeta,\mbox{\boldmath$\psi$})
=(a\zeta,\zeta)_{L^2}
 +\sum_{k=1}^n(A\partial_k\mbox{\boldmath$\psi$},\partial_k\mbox{\boldmath$\psi$})_{L^2}
 +\|\mbox{\boldmath$\psi$}'\|^2.
\]

In the following sections we use this symmetric structure of the Isobe--Kakinuma model 
to construct the solution and derive an energy estimate.

%------------------------------------------------------------------------------
\section{Construction of the solution to a reduced system}
\label{section:construction}
\setcounter{equation}{0}
\setcounter{theorem}{0}
In this section we transform the Isobe--Kakinuma model \eqref{intro:IK model} into a system of equations 
for which the hypersurface $t=0$ is noncharacteristic by using the necessary condition \eqref{intro:compatibility} 
for the existence of the solution and construct the solution of the initial value problem to the transformed 
system by using a standard parabolic regularization. 

By using the notation introduced in the previous section, the Isobe--Kakinuma model \eqref{intro:IK model} 
can be written simply as 
\begin{equation}\label{construction:IK model}
\left\{
 \begin{array}{l}
  \displaystyle
  H^{p_i}\partial_t\eta-\sum_{j=0}^NL_{ij}\phi_j=0 \qquad\mbox{for}\quad i=0,1,\ldots,N, \\
  \displaystyle
  \sum_{j=0}^NH^{p_j}\partial_t\phi_j=F_0,
 \end{array}
\right.
\end{equation}
where 
\begin{equation}\label{construction:F0}
F_0=-g\eta-\frac12(|\mbox{\boldmath$u$}|^2+w^2).
\end{equation}
The necessary condition \eqref{intro:compatibility} can also be written simply as 
\[
\sum_{j=0}^N(L_{ij}-H^{p_i}L_{0j})\phi_j=0 \qquad\mbox{for}\quad i=1,\ldots,N.
\]

We first derive an evolution equation for $\mbox{\boldmath$\phi$}=(\phi_0,\ldots,\phi_N)^{\rm T}$. 
Since the left-hand side of the above equation does not contain the term $\nabla H$, 
differentiating this with respect to the time $t$ we obtain 
\[
\sum_{j=0}^N(L_{ij}-H^{p_i}L_{0j})\partial_t\phi_j=f_i\partial_t\eta \qquad\mbox{for}\quad i=1,\ldots,N,
\]
where 
\begin{align}\label{construction:fi}
f_i &= -\frac{\partial}{\partial H}\sum_{j=0}^N(L_{ij}-H^{p_i}L_{0j})\phi_j \\
&= p_iH^{p_i-1}\biggl\{\nabla b\cdot\mbox{\boldmath$u$}-w-\sum_{j=0}^N\biggl(
 \frac{1}{p_j+1}H^{p_j+1}\Delta\phi_j-\frac{p_j}{p_j}H^{p_j}\nabla\cdot(\phi_j\nabla b)\biggr)\biggr\}
 \nonumber
\end{align}
for $i=1,\ldots,N$. 
We use the first equation in \eqref{construction:IK model} with $i=0$ to remove the time derivative 
$\partial_t\eta$ from the above equation and obtain 
\begin{equation}\label{construction:reduced 1}
\sum_{j=0}^N(L_{ij}-H^{p_i}L_{0j})\partial_t\phi_j=F_i \qquad\mbox{for}\quad i=1,\ldots,N,
\end{equation}
where 
\begin{equation}\label{consruction:Fi}
F_i=f_i\sum_{j=0}^NL_{0j}\phi_j  \qquad\mbox{for}\quad i=1,\ldots,N.
\end{equation}
By using the operator $\mathscr{L}$ introduced by \eqref{linear:L2}, the second equation in 
\eqref{construction:IK model} and \eqref{construction:reduced 1} constitute the equation 
\[
\mathscr{L}\partial_t\mbox{\boldmath$\phi$}=\mbox{\boldmath$F$},
\]
where $\mbox{\boldmath$F$}=(F_0,\ldots,F_N)^{\rm T}$. 
This is the evolution equation for $\mbox{\boldmath$\phi$}$. 
We proceed to derive an appropriate evolution equation for $\eta$. 
Let $\mbox{\boldmath$q$}=(q_0,\ldots,q_N)^{\rm T}=(q_0(H),\ldots,q_N(H))^{\rm T}$ be the 
rational functions of $H$ defined by \eqref{linear:inverse}. 
In view of the first relation in \eqref{linear:q} and the arguments in the previous section, 
we multiply the first equation in \eqref{construction:IK model} by $q_i$ and adding the 
resulting equations over $i=0,1,\ldots,N$ to obtain 
\[
\partial_t\eta=F_{N+1},
\]
where 
\begin{equation}\label{construction:FN+1}
F_{N+1}=-\sum_{i,j=0}^Nq_iL_{ij}\phi_j.
\end{equation}
To summarize we have reduced the Isobe--Kakinuma model \eqref{construction:IK model} to 
\begin{equation}\label{construction:reduced 2}
\left\{
 \begin{array}{l}
  \mathscr{L}\partial_t\mbox{\boldmath$\phi$}=\mbox{\boldmath$F$}, \\
  \partial_t\eta=F_{N+1}.
 \end{array}
\right.
\end{equation}
We note that $\mbox{\boldmath$F$}$ and $F_{N+1}$ do not contain any time derivatives and that 
$\mathscr{L}$ is invertible thanks to Lemma \ref{linear:lemma 3}. 
Therefore, the hypersurface $t=0$ is not characteristic any more for the above reduced equations. 
In the rest of this section we consider the initial value problem to \eqref{construction:reduced 2} 
under the initial conditions 
\begin{equation}\label{construction:initial conditions}
(\eta,\phi_0,\ldots,\phi_N)=(\eta_{(0)},\phi_{0(0)},\ldots,\phi_{N(0)}) \quad\makebox[3em]{at} t=0.
\end{equation}

\begin{remark}\label{construction:remark}
As mentioned in Remark {\rm \ref{intro:remark} (2)}, we have to express 
$\partial_t\phi_1(x,0),\ldots,\partial_t\phi_N(x,0)$ in terms of the initial data and $b$. 
Such an expression is given by 
$
\partial_t\mbox{\boldmath$\phi$}(\cdot,0)=\mathscr{L}^{-1}\mbox{\boldmath$F$}|_{t=0}.
$
\end{remark}

We will construct the solution to the initial value problem for the reduced system 
\eqref{construction:reduced 2}--\eqref{construction:initial conditions} by a standard 
parabolic regularization of the equations, that is, 
\begin{equation}\label{construction:regularized}
\left\{
 \begin{array}{l}
  \mathscr{L}(\partial_t\mbox{\boldmath$\phi$}-\varepsilon\Delta\mbox{\boldmath$\phi$})
   =\mbox{\boldmath$F$}, \\
  \partial_t\eta-\varepsilon\Delta\eta=F_{N+1},
 \end{array}
\right.
\end{equation}
where $\varepsilon>0$ is a small regularized parameter. 
By the definitions \eqref{construction:F0}, \eqref{consruction:Fi}, and \eqref{construction:FN+1} of 
$\mbox{\boldmath$F$}=(F_0,\ldots,F_N)^{\rm T}$ and $F_{N+1}$ (see also \eqref{construction:fi}, 
\eqref{linear:L}, and \eqref{linear:uw}) we see that 
\begin{equation}\label{construction:FF}
\left\{
 \begin{array}{l}
  F_0=F_0(\eta,H,\mbox{\boldmath$\phi$}',\nabla\mbox{\boldmath$\phi$},\nabla b), \\
  F_i=F_i(H,\nabla H,\mbox{\boldmath$\phi$}',\nabla\mbox{\boldmath$\phi$},\Delta\mbox{\boldmath$\phi$},
   \nabla b,\Delta b) \quad\mbox{for}\quad i=1,\ldots,N+1,
 \end{array}
\right.
\end{equation}
where $\mbox{\boldmath$\phi$}'=(\phi_1,\ldots,\phi_N)^{\rm T}$. 
Therefore, thanks to Lemma \ref{linear:lemma 3}, $\mathscr{L}^{-1}\mbox{\boldmath$F$}$ behaves as if 
it is a function of $(\eta,H,\mbox{\boldmath$\phi$}',\nabla\mbox{\boldmath$\phi$},\nabla b)$, 
so that we can show the following lemma.

\begin{lemma}\label{construction:existence 1}
Let $g, h, c_0, M_0$ be positive constants and $m$ an integer such that $m>n/2+1$. 
Suppose that the initial data $(\eta_{(0)},\phi_{0(0)},\ldots,\phi_{N(0)})$ and $b$ satisfy the conditions 
in \eqref{intro:conditions}, then for any $\varepsilon>0$ there exists a maximal existence time 
$T_{\varepsilon}>0$ such that the initial value problem \eqref{construction:regularized} and 
\eqref{construction:initial conditions} has a unique solution 
$(\eta^{\varepsilon},\mbox{\boldmath$\phi$}^{\varepsilon})$ satisfying 
\[
\eta^{\varepsilon},\nabla\phi_0^{\varepsilon}\in C([0,T_{\varepsilon});H^m), \quad
\phi_1^{\varepsilon},\ldots,\phi_N^{\varepsilon}\in C([0,T_{\varepsilon});H^{m+1}).
\]
\end{lemma}

We proceed to derive uniform estimates of the solution $(\eta^{\varepsilon},\mbox{\boldmath$\phi$}^{\varepsilon})$ 
with respect to the regularized parameter $\varepsilon\in(0,1]$ for a time interval $[0,T]$ independent of 
$\varepsilon$. 
To this end, we make use of a good symmetric structure of the Isobe--Kakinuma model 
discussed in the previous section. 
In order to simplify the notation we write $(\eta,\mbox{\boldmath$\phi$})$ in place of 
$(\eta^{\varepsilon},\mbox{\boldmath$\phi$}^{\varepsilon})$ in the following. 
Concerning the generalized Rayleigh--Taylor sign condition, a regularized version of the function $a$ 
defined by \eqref{intro:a} is given by 
\begin{align}\label{construction:a}
a^{\varepsilon} &= g+\sum_{i=0}^Np_iH^{p_i-1}(\partial_t\phi_i-\varepsilon\Delta\phi_i) \\
&\quad\;
+\frac12\sum_{i,j=0}^N\Bigl\{
 (p_i+p_j)H^{p_i+p_j-1}\nabla\phi_i\cdot\nabla\phi_j
 -2p_i(p_i+p_j-1)H^{p_i+p_j-2}\phi_i\nabla b\cdot\nabla\phi_j \nonumber \\
&\makebox[6em]{}
+p_ip_j(p_i+p_j-2)H^{p_i+p_j-3}(1+|\nabla b|^2)\phi_i\phi_j\Bigr\}. \nonumber
\end{align}
Here, we have $\partial_t\mbox{\boldmath$\phi$}-\varepsilon\Delta\mbox{\boldmath$\phi$}
=\mathscr{L}^{-1}\mbox{\boldmath$F$}$ so that 
$a^{\varepsilon}(x,0)=a(x,0)$. See also Remark \ref{construction:remark}. 
Therefore, by the sign condition in \eqref{intro:conditions} we have 
\begin{equation}\label{construction:sign}
a^{\varepsilon}(x,0)\geq c_0 \quad\mbox{for}\quad x\in\mathbf{R}^n.
\end{equation}

We will derive a regularized version of the symmetric positive system \eqref{linear:positive system} 
for $(\eta,\mbox{\boldmath$\phi$})$. 
First, we derive an evolution equation for $\mbox{\boldmath$\phi$}$. 
The first equation in \eqref{construction:regularized} is written as 
\begin{equation}\label{construction:regularized 2}
\left\{
 \begin{array}{l}
  \mbox{\boldmath$l$}\cdot(\partial_t\mbox{\boldmath$\phi$}-\varepsilon\Delta\mbox{\boldmath$\phi$})=F_0, \\[0.5ex]
  \mathscr{L}_i(\partial_t\mbox{\boldmath$\phi$}-\varepsilon\Delta\mbox{\boldmath$\phi$})=F_i
   \quad\mbox{for}\quad i=1,\ldots,N,
 \end{array}
\right.
\end{equation}
where $\mbox{\boldmath$l$}=\mbox{\boldmath$l$}(H)=(H^{p_0},\ldots,H^{p_N})^{\rm T}$. 
Here we remark that $F_i$ $(i=1,\ldots,N)$ is a correction of lower order terms and its explicit form has 
no importance whereas $F_0$ contains principal terms and we have to treat it carefully. 
In other words, the second equation in \eqref{construction:regularized 2} is quasilinear whereas 
the first one is fully nonlinear. 
Applying $\Delta$ to the first equation in \eqref{construction:regularized 2} we have 
\begin{align*}
& \mbox{\boldmath$l$}(H)\cdot\Delta
 (\partial_t\mbox{\boldmath$\phi$}-\varepsilon\Delta\mbox{\boldmath$\phi$})
+(\mbox{\boldmath$l$}^{(1)}(H)\cdot(\partial_t\mbox{\boldmath$\phi$}-\varepsilon\Delta\mbox{\boldmath$\phi$}))
 \Delta H-\Delta F_0 \\
& = -([\Delta,\mbox{\boldmath$l$}(H)]-\mbox{\boldmath$l$}^{(1)}(H)(\Delta H))\cdot
 \mathscr{L}^{-1}\mbox{\boldmath$F$},
\end{align*}
where $\mbox{\boldmath$l$}^{(1)}(H)=(0,p_1H^{p_1-1},\ldots,p_NH^{p_N-1})^{\rm T}$ is the derivative of 
$\mbox{\boldmath$l$}(H)$ with respect to $H$. 
Here we see that 
\begin{align*}
-\Delta F_0
&= g\Delta\eta+\mbox{\boldmath$u$}\cdot\Delta\mbox{\boldmath$u$}+w\Delta w
 +\frac12\bigl\{(\Delta(\mbox{\boldmath$u$}\cdot\mbox{\boldmath$u$})
  -2\mbox{\boldmath$u$}\cdot\Delta\mbox{\boldmath$u$})
  +(\Delta(w^2)-2w\Delta w)\bigr\},
\end{align*}
where 
\begin{align*}
\Delta\mbox{\boldmath$u$}
&= \sum_{i=0}^N\bigl\{  H^{p_i}\nabla\Delta\phi_i
 +p_i(H^{p_i-1}\nabla\phi_i-(p_i-1)H^{p_i-2}\phi_i\nabla b)(\Delta\eta) \\
&\qquad\quad
 +\bigl([\Delta,H^{p_i}]-p_iH^{p_i-1}(\Delta H)-p_iH^{p_i-1}
  (\Delta b)\bigr)\nabla\phi_i
  -p_iH^{p_i-1}\Delta(\phi_i\nabla b) \\
&\qquad\quad
 -p_i\bigl([\Delta,H^{p_i-1}]-(p_i-1)H^{p_i-2}(\Delta H)-(p_i-1)H^{p_i-2}
  (\Delta b)\bigr)(\phi_i\nabla b)\bigr\},
\end{align*}
\begin{align*}
\Delta w
&= \sum_{i=0}^N\bigl\{ p_i(p_i-1)H^{p_i-2}\phi_i\Delta\eta
 +p_iH^{p_i-1}\Delta\phi_i \\
&\qquad\quad
 +p_i\bigl([\Delta,H^{p_i-1}]-(p_i-1)H^{p_i-2}(\Delta H)-(p_i-1)H^{p_i-2}
  (\Delta b)\bigr)\phi_i \bigr\}.
\end{align*}
Therefore, we obtain 
\begin{equation}\label{construction:quasilinear 1}
\mbox{\boldmath$l$}(H)\cdot\Delta
 (\partial_t\mbox{\boldmath$\phi$}-\varepsilon\Delta\mbox{\boldmath$\phi$})
+ \mbox{\boldmath$l$}(H)\cdot(\mbox{\boldmath$u$}\cdot\nabla)\Delta\mbox{\boldmath$\phi$}
+ a^{\varepsilon}\Delta\eta
 = G_{00},
\end{equation}
where 
\begin{align}\label{construction:G00}
G_{00}
&= -([\Delta,\mbox{\boldmath$l$}(H)]-\mbox{\boldmath$l$}^{(1)}(H)(\Delta H))\cdot
  \mathscr{L}^{-1}\mbox{\boldmath$F$}
 +(\mbox{\boldmath$l$}^{(1)}(H)\cdot(\mathscr{L}^{-1}\mbox{\boldmath$F$}))
 \Delta b \\
&\quad\;
 -\frac12\bigl\{(\Delta(\mbox{\boldmath$u$}\cdot\mbox{\boldmath$u$})
  -2\mbox{\boldmath$u$}\cdot\Delta\mbox{\boldmath$u$})
  +(\Delta(w^2)-2w\Delta w)\bigr\} \nonumber \\
&\quad\;
 -\mbox{\boldmath$u$}\cdot\sum_{i=0}^N\bigl\{ 
 \bigl([\Delta,H^{p_i}]-p_iH^{p_i-1}(\Delta H)-p_iH^{p_i-1}
  (\Delta b)\bigr)\nabla\phi_i
  -p_iH^{p_i-1}\Delta(\phi_i\nabla b) \nonumber \\
&\qquad\qquad\quad
 -p_i\bigl([\Delta,H^{p_i-1}]-(p_i-1)H^{p_i-2}(\Delta H)-(p_i-1)H^{p_i-2}
  (\Delta b)\bigr)(\phi_i\nabla b)\bigr\} \nonumber \\
&\quad\;
 -w\sum_{i=0}^N\bigl\{ p_iH^{p_i-1}\Delta\phi_i \nonumber \\
&\qquad\qquad\quad
 +p_i\bigl([\Delta,H^{p_i-1}]-(p_i-1)H^{p_i-2}(\Delta H)-(p_i-1)H^{p_i-2}
  (\Delta b)\bigr)\phi_i \bigr\}. \nonumber
\end{align}
In view of \eqref{linear:L2} and \eqref{linear:L} we divide the operator $\mathscr{L}_i$ into 
its principal part and the remainder part $\mathscr{L}_i^{\rm low}$ as 
\[
\mathscr{L}_i\mbox{\boldmath$\psi$}
= (-\mbox{\boldmath$a$}_{i}(H)+H^{p_i}\mbox{\boldmath$a$}_{0}(H))\cdot\Delta\mbox{\boldmath$\psi$}
 +\mathscr{L}_i^{\rm low}\mbox{\boldmath$\psi$} \quad\mbox{for}\quad i=1,\ldots,N,
\]
where $A(H)=(\mbox{\boldmath$a$}_{0}(H),\ldots,\mbox{\boldmath$a$}_{N}(H))$ is the matrix defined by 
\eqref{linear:A} and 
\begin{align}\label{construction:Lilow}
\mathscr{L}_i^{\rm low}\mbox{\boldmath$\psi$}
&= \sum_{j=0}^N\biggl(-\frac{p_ip_j}{p_j(p_i+p_j)}H^{p_i+p_j}\nabla\cdot(\psi_j\nabla b)
 -\frac{p_i}{p_i+p_j}H^{p_i+p_j}\nabla b\cdot\nabla\psi_j \\
&\phantom{ = \sum\biggl( }
 +\frac{p_ip_j}{p_i+p_j-1}H^{p_i+p_j-1}(1+|\nabla b|^2)\psi_j\biggr). \nonumber
\end{align}
Therefore, we obtain 
\begin{equation}\label{construction:quasilinear 2}
(\mbox{\boldmath$a$}_{i}(H)-H^{p_i}\mbox{\boldmath$a$}_{0}(H))\cdot\Delta
 (\partial_t\mbox{\boldmath$\phi$}-\varepsilon\Delta\mbox{\boldmath$\phi$})
= G_{0i} \quad\mbox{for}\quad i=1,\ldots,N,
\end{equation}
where 
\begin{align}\label{construction:G0i}
G_{0i}
&= -F_i
 +\mathscr{L}_i^{\rm low}\mathscr{L}^{-1}\mbox{\boldmath$F$}.
\end{align}
We can rewrite \eqref{construction:quasilinear 1} and \eqref{construction:quasilinear 2} 
in a matrix form as 
\[
\tilde{A}(H)
\left(
 \begin{array}{c}
  \mbox{\boldmath$a$}_0\cdot\Delta(\partial_t\mbox{\boldmath$\phi$}-\varepsilon\Delta\mbox{\boldmath$\phi$}) \\
  \Delta(\partial_t\mbox{\boldmath$\phi$}-\varepsilon\Delta\mbox{\boldmath$\phi$})
 \end{array}
\right)
+ 
\left(
 \begin{array}{c}
  \mbox{\boldmath$l$}(H)^{\rm T}(\mbox{\boldmath$u$}\cdot\nabla)\Delta\mbox{\boldmath$\phi$}
   +a^{\varepsilon}\Delta\eta \\
  \mbox{\boldmath$0$}
 \end{array}
\right) 
= 
\left(
 \begin{array}{c}
  G_{00} \\
  \mbox{\boldmath$G$}_0
 \end{array}
\right),
\]
where $\tilde{A}(H)$ is the matrix in \eqref{linear:tildeA} and 
$\mbox{\boldmath$G$}_0=(0,G_{01},\ldots,G_{0N})^{\rm T}$. 
Therefore, using the notation in \eqref{linear:inverse} and the relations in \eqref{linear:q} we have 
\[
A(H)\Delta(\partial_t\mbox{\boldmath$\phi$}-\varepsilon\Delta\mbox{\boldmath$\phi$})
+ q(H)\mbox{\boldmath$l$}(H)\{ \mbox{\boldmath$l$}(H)^{\rm T}(\mbox{\boldmath$u$}\cdot\nabla)
   \Delta\mbox{\boldmath$\phi$} +a^{\varepsilon}\Delta\eta\} = \mbox{\boldmath$G$},
\]
where 
\begin{equation}\label{construction:G}
\mbox{\boldmath$G$}
= q(H)\mbox{\boldmath$l$}(H)G_{00}+A(H)Q(H)\mbox{\boldmath$G$}_0.
\end{equation}
This is the desired equation for $\mbox{\boldmath$\phi$}$. 

Secondly, we derive an evolution equation for $\eta$. 
In view of \eqref{linear:L} we divide the operator $L_{ij}$ into its principal part and 
the remainder part $L_{ij}^{\rm low}$ as 
\[
L_{ij}\phi_j=-a_{ij}(H)\Delta\phi_j-H^{p_i}(H^{p_j}\nabla\phi_j-p_jH^{p_j-1}\phi_j\nabla b)\cdot\nabla\eta
 +L_{ij}^{\rm low}\phi_j,
\]
where 
\begin{align}\label{construction:Llow}
L_{ij}^{\rm low}\phi_j
&= H^{p_i}(H^{p_j}\nabla\phi_j-p_jH^{p_j-1}\phi_j\nabla b)\cdot\nabla b
 + \frac{p_j}{p_i+p_j}H^{p_i+p_j}\nabla\cdot(\phi_j\nabla b) \\
&\quad\;
 -\frac{p_i}{p_i+p_j}H^{p_i+p_j}\nabla b\cdot\nabla\phi_j
 +\frac{p_ip_j}{p_i+p_j-1}H^{p_i+p_j-1}(1+|\nabla b|^2)\phi_j. \nonumber
\end{align}
Thanks to the relations in \eqref{linear:q}, the second equation in \eqref{construction:regularized} 
can be written as 
\[
\partial_t\eta-\varepsilon\Delta\eta
 +\mbox{\boldmath$u$}\cdot\nabla\eta
 +q(H)\mbox{\boldmath$l$}(H)^{\rm T}\Delta\mbox{\boldmath$\phi$}=G_0,
\]
where 
\begin{equation}\label{construction:G0}
G_0 = -\sum_{i,j=0}^Nq_i(H)L_{ij}^{\rm low}\phi_j.
\end{equation}
This is the desired equation for $\eta$.

To summarize, we have derived the equations 
\begin{equation}\label{construction:quasilinear 3}
\left\{
 \begin{array}{l}
  \partial_t\eta-\varepsilon\Delta\eta
   +\mbox{\boldmath$u$}\cdot\nabla\eta
   +q(H)\mbox{\boldmath$l$}(H)^{\rm T}\Delta\mbox{\boldmath$\phi$}=G_0, \\[0.5ex]
  A(H)\Delta(\partial_t\mbox{\boldmath$\phi$}-\varepsilon\Delta\mbox{\boldmath$\phi$})
   + q(H)\mbox{\boldmath$l$}(H)\{ \mbox{\boldmath$l$}(H)^{\rm T}(\mbox{\boldmath$u$}\cdot\nabla)
    \Delta\mbox{\boldmath$\phi$} +a^{\varepsilon}\Delta\eta\} = \mbox{\boldmath$G$}.
 \end{array}
\right.
\end{equation}
Using this we can derive uniform estimate of the solution 
$(\eta^{\varepsilon},\mbox{\boldmath$\phi$}^{\varepsilon})$ of the initial value problem 
to the regularized equation \eqref{construction:regularized}.

\begin{lemma}\label{construction:uniform}
Under the hypothesis of Lemma {\rm \ref{construction:existence 1}}, 
there exist a time $T>0$ and a constant $C>0$ independent of $\varepsilon$ such that the solution 
$(\eta^{\varepsilon},\mbox{\boldmath$\phi$}^{\varepsilon})$ obtained in Lemma {\rm \ref{construction:existence 1}}
satisfies the uniform estimate 
\begin{align}\label{construction:uniform estimate}
& \sup_{0\leq t\leq T}\bigl(
\|\eta^{\varepsilon}(t)\|_m^2+\|\nabla\phi_0^{\varepsilon}(t)\|_m^2
 +\|(\phi_1^{\varepsilon}(t),\ldots,\phi_N^{\varepsilon}(t))\|_{m+1}^2\bigr) \\
&\quad
 +\varepsilon\int_0^T\bigl(\|\eta^{\varepsilon}(t)\|_{m+1}^2+\|\nabla\phi_0^{\varepsilon}(t)\|_{m+1}^2
 +\|(\phi_1^{\varepsilon}(t),\ldots,\phi_N^{\varepsilon}(t))\|_{m+2}^2\bigr){\rm d}t \leq C
 \nonumber
\end{align}
for $0<\varepsilon\leq 1$. 
\end{lemma}

\noindent
{\bf Proof}. \ 
Once again we simply write $(\eta,\mbox{\boldmath$\phi$})$ in place of 
$(\eta^{\varepsilon},\mbox{\boldmath$\phi$}^{\varepsilon})$. 
Let $\alpha=(\alpha_1,\ldots,\alpha_n)$ be a multi-index satisfying $1\leq |\alpha|\leq m$. 
Without loss of generality, we can assume that $\alpha_1\geq1$. 
Applying the differential operator $a^{\varepsilon}\partial^{\alpha}$ to the first equation in 
\eqref{construction:quasilinear 3} we have 
\begin{equation}\label{construction:quasilinear 4}
a^{\varepsilon}\partial_t\partial^{\alpha}\eta-\varepsilon\nabla\cdot(a^{\varepsilon}\nabla\partial^{\alpha}\eta)
 +a^{\varepsilon}\mbox{\boldmath$u$}\cdot\nabla\partial^{\alpha}\eta
 +\sum_{k=1}^n\partial_k(a^{\varepsilon}q\mbox{\boldmath$l$}^{\rm T}\partial_k\partial^{\alpha}
  \mbox{\boldmath$\phi$})=F_{0,\alpha},
\end{equation}
where 
\begin{equation}\label{construction:F00}
F_{0,\alpha} = a^{\varepsilon}(\partial^{\alpha}G_0-[\partial^{\alpha},\mbox{\boldmath$u$}]\cdot\nabla\eta
 -[\partial^{\alpha},q\mbox{\boldmath$l$}^{\rm T}]\Delta\mbox{\boldmath$\phi$})
 -\varepsilon\nabla a^{\varepsilon}\cdot\nabla\partial^{\alpha}\eta
 +\sum_{k=1}^n(\partial_k(a^{\varepsilon}q\mbox{\boldmath$l$}^{\rm T}))
  \partial_k\partial^{\alpha}\mbox{\boldmath$\phi$}.
\end{equation}
Here and in what follows, for simplicity, we omit the dependence of $H$ in the notation. 
Applying the differential operator $\partial^{\alpha}$ to the second equation in 
\eqref{construction:quasilinear 3} we have 
\begin{align}\label{construction:quasilinear 5}
-\sum_{k=1}^n\partial_k(A\partial_k\partial_t\partial^{\alpha}\mbox{\boldmath$\phi$})
 +\varepsilon\Delta(A\Delta\partial^{\alpha}\mbox{\boldmath$\phi$})
 -& \sum_{k=1}^n\partial_k\bigl\{q\mbox{\boldmath$l$}\mbox{\boldmath$l$}^{\rm T}(\mbox{\boldmath$u$}\cdot\nabla)
  \partial_k\partial^{\alpha}\mbox{\boldmath$\phi$}
 +a^{\varepsilon}q\mbox{\boldmath$l$}\partial_k\partial^{\alpha}\eta\bigr\} \\
&= \sum_{k=1}^n\partial_k\mbox{\boldmath$F$}_{k,\alpha}, \nonumber
\end{align}
where 
\begin{align}\label{construction:F}
\mbox{\boldmath$F$}_{k,\alpha}
&= -\delta_{1k}\partial^{\alpha'}\biggl\{ \mbox{\boldmath$G$}
   +\sum_{l=1}^n\bigl((\partial_lA)\partial_l\mathscr{L}^{-1}\mbox{\boldmath$F$}
    +[\partial_l,q\mbox{\boldmath$l$}\mbox{\boldmath$l$}^{\rm T}(\mbox{\boldmath$u$}\cdot\nabla)]
     \partial_l\mbox{\boldmath$\phi$}
    +(\partial_l(a^{\varepsilon}q\mbox{\boldmath$l$}))\partial_l\eta\bigr)\biggr\} \\
&\quad\;
 +\varepsilon(\partial_kA)\Delta\partial^{\alpha}\mbox{\boldmath$\phi$}
 +[\partial^{\alpha},A]\partial_k\mathscr{L}^{-1}\mbox{\boldmath$F$}
 +[\partial^{\alpha},q\mbox{\boldmath$l$}\mbox{\boldmath$l$}^{\rm T}(\mbox{\boldmath$u$}\cdot\nabla)]
  \partial_k\mbox{\boldmath$\phi$}
 +[\partial^{\alpha},a^{\varepsilon}q\mbox{\boldmath$l$}]\partial_k\eta.
    \nonumber
\end{align}
Here, $\alpha'=(\alpha_1-1,\alpha_2,\ldots,\alpha_n)$ and $\delta_{1k}$ is the Kronecker delta. 
\eqref{construction:quasilinear 4} and \eqref{construction:quasilinear 5} can be written in the 
matrix form as 
\begin{align}\label{construction:positive system 2}
& \left(
 \begin{array}{cc}
  a^{\varepsilon} & \mbox{\boldmath$0$}^{\rm T} \\
  \mbox{\boldmath$0$} & -\sum_{k=1}^n\partial_k(A\partial_k\,\cdot\,)
 \end{array}
\right)
\partial_t
\left(
 \begin{array}{c}
  \partial^{\alpha}\eta \\
  \partial^{\alpha}\mbox{\boldmath$\phi$}
 \end{array}
\right)
+\varepsilon
\left(
 \begin{array}{cc}
  -\nabla\cdot(a^{\varepsilon}\nabla\,\cdot\,) & \mbox{\boldmath$0$}^{\rm T} \\
  \mbox{\boldmath$0$} & \Delta(A\Delta\,\cdot\,)
 \end{array}
\right)
\left(
 \begin{array}{c}
  \partial^{\alpha}\eta \\
  \partial^{\alpha}\mbox{\boldmath$\phi$}
 \end{array}
\right) \\
&\quad +
\left(
 \begin{array}{cc}
  a^{\varepsilon}\mbox{\boldmath$u$}\cdot\nabla
   & \sum_{k=1}^n\partial_k(a^{\varepsilon}q\mbox{\boldmath$l$}^{\rm T}\partial_k\,\cdot\,) \\
  -\sum_{k=1}^n\partial_k(a^{\varepsilon}q\mbox{\boldmath$l$}\partial_k\,\cdot\,)
   & -\sum_{k=1}^n\partial_k(q\mbox{\boldmath$l$}\mbox{\boldmath$l$}^{\rm T}(\mbox{\boldmath$u$}\cdot\nabla)
  \partial_k\,\cdot\,)
 \end{array}
\right)
\left(
 \begin{array}{c}
  \partial^{\alpha}\eta \\
  \partial^{\alpha}\mbox{\boldmath$\phi$}
 \end{array}
\right) \nonumber \\
&\quad =
\left(
 \begin{array}{c}
  F_{0,\alpha} \\
  \sum_{k=1}^n\partial_k\mbox{\boldmath$F$}_{k,\alpha}
 \end{array}
\right), \nonumber
\end{align}
which forms a symmetric positive system.

In view of \eqref{construction:positive system 2} we define an energy function $\mathscr{E}_m(t)$ and a 
dissipation function $\mathscr{D}_m(t)$ by 
\begin{align*}
& \mathscr{E}_m(t) = \sum_{|\alpha|\leq m}\biggl\{
 (a^{\varepsilon}\partial^{\alpha}\eta(t),\partial^{\alpha}\eta(t))_{L^2}
 +\sum_{k=1}^n(A\partial_k\partial^{\alpha}\mbox{\boldmath$\phi$}(t),
  \partial_k\partial^{\alpha}\mbox{\boldmath$\phi$}(t))_{L^2}\biggr\}
 +\|\mbox{\boldmath$\phi$}'(t)\|^2, \\
& \mathscr{D}_m(t) = \sum_{|\alpha|\leq m}\bigl\{
 (a^{\varepsilon}\nabla\partial^{\alpha}\eta(t),\nabla\partial^{\alpha}\eta(t))_{L^2}
 +(A\Delta\partial^{\alpha}\mbox{\boldmath$\phi$}(t),\Delta\partial^{\alpha}\mbox{\boldmath$\phi$}(t))_{L^2}
 \bigr\},
\end{align*}
which will be equivalent to 
\[
E_m(t) = \|\eta(t)\|_m^2+\|\nabla\phi_0(t)\|_m^2+\|\mbox{\boldmath$\phi$}'(t)\|_{m+1}^2, \quad
D_m(t) = \|\nabla\eta(t)\|_m^2+\|\Delta\mbox{\boldmath$\phi$}(t)\|_m^2,
\]
respectively, 
where $\mbox{\boldmath$\phi$}'=(\phi_1,\ldots,\phi_N)^{\rm T}$. 
In view of \eqref{intro:conditions}, \eqref{construction:FF}, \eqref{construction:a}, \eqref{construction:sign}, 
and Lemma \ref{linear:lemma 3}, we see that there exists a constant $C_0=C_0(g,h,c_0,M_0)>0$ such that 
\[
c_0\leq H(x,0)\leq C_0, \quad c_0\leq a^{\varepsilon}(x,0)\leq C_0 
 \qquad\mbox{for}\quad x\in\mathbf{R}^n.
\]
Now, we assume that 
\begin{equation}\label{construction:assumption}
E_m(t)+\varepsilon\int_0^tE_{m+1}(\tau){\rm d}\tau \leq M_1, \quad
 \frac{c_0}{2}\leq H(x,t)\leq 2C_0, \quad 
 \frac{c_0}{2}\leq a^{\varepsilon}(x,t)\leq 2C_0 
\end{equation}
for $0\leq t\leq T$ and $x\in\mathbf{R}^n$, where the constant $M_1$ and the time $T$ will be determined later. 
In the following we simply write the constants depending only on $(g,h,c_0,C_0,M_0)$ by $C_1$ and the 
constants depending also on $M_1$ by $C_2$, which may change from line to line. 
Then, it holds that 
\[
C_1^{-1}E_m(t) \leq \mathscr{E}_m(t) \leq C_1E_m(t), \quad 
C_1^{-1}D_m(t) \leq \mathscr{D}_m(t) \leq C_1D_m(t)
\]
for $0\leq t\leq T$. 
We are going to evaluate the evolution of the energy function $\mathscr{E}_m(t)$. 
To this end, we make use of the symmetric form \eqref{construction:positive system 2} for the case 
$1\leq|\alpha| \leq m$ and of \eqref{construction:regularized} directly for the case $|\alpha|=0$. 
Then, by integration by parts we see that 
\begin{align}\label{construction:energy}
&\frac{\rm d}{{\rm d}t}\mathscr{E}_m(t)+2\varepsilon(\mathscr{D}_m(t)+\|\nabla\mbox{\boldmath$\phi$}'(t)\|^2) \\
&= \sum_{|\alpha|\leq m}\biggl\{
  ((\partial_ta^{\varepsilon})\partial^{\alpha}\eta,\partial^{\alpha}\eta)_{L^2}
 +\sum_{k=1}^n((\partial_tA)\partial_k\partial^{\alpha}\mbox{\boldmath$\phi$},
  \partial_k\partial^{\alpha}\mbox{\boldmath$\phi$})_{L^2}\biggr\} \nonumber\\
&\quad\;
 +\sum_{1\leq |\alpha|\leq m}\biggr\{
  ((\nabla\cdot(a^{\varepsilon}\mbox{\boldmath$u$}))\partial^{\alpha}\eta,\partial^{\alpha}\eta)_{L^2}
 +\sum_{k=1}^n(\biggl(\sum_{l=1}^n\partial_l(u_lq\mbox{\boldmath$l$}\mbox{\boldmath$l$}^{\rm T})\biggr)
  \partial_k\partial^{\alpha}\mbox{\boldmath$\phi$},
  \partial_k\partial^{\alpha}\mbox{\boldmath$\phi$})_{L^2} \nonumber\\
&\phantom{\quad\; +\sum_{1\leq |\alpha|\leq m}\biggr\{ }
 +2(F_{0,\alpha},\partial^{\alpha}\eta)_{L^2}
 -2\sum_{k=1}^n(\mbox{\boldmath$F$}_{k,\alpha},\partial_k\partial^{\alpha}\mbox{\boldmath$\phi$})_{L^2}
 \biggr\} \nonumber\\
&\quad\;
 +2(a^{\varepsilon}F_{N+1},\eta)_{L^2}
 +2\sum_{k=1}^n(\partial_k\mathscr{L}^{-1}\mbox{\boldmath$F$},A\partial_k\mbox{\boldmath$\phi$})_{L^2}
 +2((\mathscr{L}^{-1}\mbox{\boldmath$F$})',\mbox{\boldmath$\phi$}')_{L^2} \nonumber\\
&\quad\;
 -2\varepsilon(\nabla\eta,\eta\nabla a^{\varepsilon})_{L^2}
 -2\varepsilon\sum_{k=1}^n(\Delta\mbox{\boldmath$\phi$},(\partial_kA)\partial_k\mbox{\boldmath$\phi$})_{L^2},
 \nonumber
\end{align}
where $(\mathscr{L}^{-1}\mbox{\boldmath$F$})'$ is the last $N$ components of 
$\mathscr{L}^{-1}\mbox{\boldmath$F$}$. 
Here, we see easily that 
$\|(\mbox{\boldmath$u$},w)\|_m\leq C_2$, $\|\mathscr{L}^{-1}\mbox{\boldmath$F$}\|_m\leq C_2$, 
and $\|\nabla a^{\varepsilon}\|_{m-1}\leq C_2$. 
It follows from \eqref{construction:F00} and \eqref{construction:F} 
(see also \eqref{construction:G00}, \eqref{construction:Lilow}, and 
\eqref{construction:G0i}--\eqref{construction:G0}) that 
$\|F_{0,\alpha}\|\leq C_2(1+\varepsilon\|\nabla\eta\|_m)$ and 
$\|\mbox{\boldmath$F$}_{k,\alpha}\|\leq C_2(1+\varepsilon\|\Delta\mbox{\boldmath$\phi$}\|_m)$ 
for $k=1,\ldots,n$. 
We need to evaluate $\partial_ta^{\varepsilon}$. 
Note that we can write formally as 
$[\partial_t,\mathscr{L}]=(\partial_t\eta)(\frac{\partial}{\partial H}\mathscr{L})$ 
because the coefficients of the differential operator $\mathscr{L}$ do not depend on $\nabla H$ but on $H$. 
By the relation 
\[
\partial_t(\partial_t\mbox{\boldmath$\phi$}-\varepsilon\Delta\mbox{\boldmath$\phi$})
= -\mathscr{L}^{-1}[\partial_t,\mathscr{L}]
 \mathscr{L}^{-1}\mbox{\boldmath$F$}+\mathscr{L}^{-1}\partial_t\mbox{\boldmath$F$}
\]
and Lemma \ref{linear:lemma 3} we see that 
\begin{align*}
\|\partial_t(\partial_t\mbox{\boldmath$\phi$}-\varepsilon\Delta\mbox{\boldmath$\phi$})\|_{m-1}
&\leq C_1\biggl(\|[\partial_t,\mathscr{L}_0]\mathscr{L}^{-1}\mbox{\boldmath$F$}\|_{m-1}
 +\sum_{i=1}^N\|[\partial_t,\mathscr{L}_i]\mathscr{L}^{-1}\mbox{\boldmath$F$}\|_{m-3} \\
&\qquad\quad
 +\|\partial_tF_0\|_{m-1}+\|(\partial_tF_1,\ldots,\partial_tF_N)\|_{m-3}\biggr) \\
&\leq C_1(\|\partial_t\mbox{\boldmath$\phi$}\|_m+\|\partial_t\eta\|_{m-1}).
\end{align*}
Therefore, we have $\|\partial_ta^{\varepsilon}\|_{m-1}\leq 
C_2(\|\partial_t\mbox{\boldmath$\phi$}\|_m+\|\partial_t\eta\|_{m-1})$. 
In view of 
$\|\partial_t\mbox{\boldmath$\phi$}-\varepsilon\Delta\mbox{\boldmath$\phi$}\|_m\leq C_2$ and 
$\|\partial_t\eta-\varepsilon\Delta\eta\|_{m-1}\leq C_2$, we also have 
$\|\partial_t\mbox{\boldmath$\phi$}\|_m\leq C_2(1+\varepsilon\|\Delta\mbox{\boldmath$\phi$}\|_m)$ and 
$\|\partial_t\eta\|_{m-1} \leq C_2(1+\varepsilon\|\nabla\eta\|_m)$. 
Thus, by \eqref{construction:energy} we have 
$\frac{\rm d}{{\rm d}t}\mathscr{E}_m(t)+2\varepsilon\mathscr{D}_m(t)
\leq C_2(1+\varepsilon(\|\Delta\mbox{\boldmath$\phi$}\|_m+\varepsilon\|\nabla\eta\|_m))$, 
which implies that 
$\mathscr{E}_m(t)+2\varepsilon\int_0^t\mathscr{D}_m(\tau){\rm d}\tau \leq C_1+C_2(t+\sqrt{t})$. 
To summarize, we have obtained the estimate 
\[
\left\{
 \begin{array}{l}
  \displaystyle
  E_m(t)+\varepsilon\int_0^tE_{m+1}(\tau){\rm d}\tau \leq C_1+C_2(t+\sqrt{t}), \\
  |H(x,t)-H(x,0)|+|a^{\varepsilon}(x,t)-a^{\varepsilon}(x,0)|\leq C_2(t+\sqrt{t})
   \quad\mbox{for}\quad x\in\mathbf{R}^n. 
 \end{array}
\right.
\]
Now, we define the constant $M_1$ by $M_1=2C_1$, and then the time $T$ sufficiently small 
so that $C_2(T+\sqrt{T})\ll1$. 
Then, we see that \eqref{construction:assumption} holds. 
The proof is complete. 
\quad$\Box$

\bigskip
Once we obtain this kind of uniform estimate \eqref{construction:uniform estimate}, 
we can pass to the limit $\varepsilon\to+0$ in the regularized problem \eqref{construction:regularized} 
and \eqref{construction:initial conditions} and obtain the following lemma.

\begin{lemma}\label{construction:existence}
Under the hypothesis of Theorem {\rm \ref{intro:theorem}}, there exists a time $T>0$ such that the initial 
value problem \eqref{construction:reduced 2}--\eqref{construction:initial conditions} has a unique 
solution $(\eta,\mbox{\boldmath$\phi$})$ satisfying 
\[
\eta,\nabla\phi_0\in C([0,T];H^m), \quad \phi_1,\ldots,\phi_N\in C([0,T];H^{m+1}).
\]
\end{lemma}

%------------------------------------------------------------------------------
\section{Proof of the main theorem}
\label{section:proof}
\setcounter{equation}{0}
\setcounter{theorem}{0}
We will show that the solution to the transformed problem 
\eqref{construction:reduced 2}--\eqref{construction:initial conditions} is the solution of the Isobe--Kakinuma 
model \eqref{intro:IK model}--\eqref{intro:initial conditions} if the initial data 
$(\eta_{(0)},\phi_{0(0)},\ldots,\phi_{N(0)})$ and the bottom topography $b$ satisfy the necessary condition 
\eqref{intro:compatibility}. 
Let $(\eta,\mbox{\boldmath$\phi$})$ be the solution of 
\eqref{construction:reduced 2}--\eqref{construction:initial conditions} obtained in Lemma \ref{construction:existence}. 
Then, we have $\mathscr{L}_0\partial_t\mbox{\boldmath$\phi$}=F_0$, which is exactly the second equation in 
\eqref{construction:IK model}. 
Therefore, it is sufficient to show that the first equation in \eqref{construction:IK model} holds for 
$i=0,1,\ldots,N$. 
To this end, putting 
\begin{equation}\label{proof:R}
R_i=H^{p_i}\partial_t\eta-\sum_{j=0}^NL_{ij}\phi_j \quad\mbox{for}\quad i=0,1,\ldots,N,
\end{equation}
we are going to show $R_i(x,t)\equiv0$ for $i=0,1,\ldots,N$. 
We also introduce auxiliary functions 
\begin{equation}\label{proof:RR}
\tilde{R}_i=\mathscr{L}_i\mbox{\boldmath$\phi$} \quad\mbox{for}\quad i=1,\ldots,N.
\end{equation}
Then, we see that the necessary condition \eqref{intro:compatibility} is equivalent to $\tilde{R}_i=0$ for 
$i=1,\ldots,N$. 
We have also the relations 
\begin{equation}\label{proof:relation 1}
\tilde{R}_i=H^{p_i}R_0-R_i  \quad\mbox{for}\quad i=1,\ldots,N.
\end{equation}
Differentiate \eqref{proof:RR} with respect to the time $t$ and using \eqref{proof:R} with $i=0$ to eliminate 
the time derivative $\partial_t\eta$ and 
the equations $\mathscr{L}_i\partial_t\mbox{\boldmath$\phi$}=F_i$, we obtain 
\begin{equation}\label{proof:ODE1}
\partial_t\tilde{R}_i=-f_iR_0 \quad\mbox{for}\quad i=1,\ldots,N,
\end{equation}
where $f_i$ is the function defined by \eqref{construction:fi}. 
Let $\mbox{\boldmath$q$}=\mbox{\boldmath$q$}(H)=(q_0(H),\ldots,q_N(H))$ be rational functions of $H$ defined 
by \eqref{linear:inverse}. 
Multiplying \eqref{proof:R} by $q_i$, adding the resulting equations over $i=0,1,\ldots,N$, and using the 
relation \eqref{linear:q} we obtain 
\[
\sum_{i=0}^Nq_iR_i=-\partial_t\eta-\sum_{i,j=0}^Nq_iL_{ij}\phi_j=0,
\]
where we used the second equation in \eqref{construction:reduced 2}. 
This together with \eqref{proof:relation 1} implies 
\begin{equation}\label{proof:relation 2}
R_0=-\sum_{j=1}^Nq_j\tilde{R}_j.
\end{equation}
Here, we used the relation $\mbox{\boldmath$l$}\cdot\mbox{\boldmath$q$}=-1$ once again. 
Plugging this into \eqref{proof:ODE1} we obtain a system of linear homogeneous ordinary differential 
equations for $(\tilde{R}_1,\ldots,\tilde{R}_N)$: 
\[
\partial_t\tilde{R}_i=f_i\sum_{j=1}^Nq_j\tilde{R}_j \quad\mbox{for}\quad i=1,\ldots,N.
\]
The necessary condition \eqref{intro:compatibility} for the initial data 
$(\eta_{(0)},\phi_{0(0)},\ldots,\phi_{N(0)})$ and the bottom topography $b$ is equivalent to 
$\tilde{R}_i(x,0)\equiv0$ for $i=1,\ldots,N$, so that the uniqueness of the solution to the 
initial value problem implies $\tilde{R}_i=0$ for $i=1,\ldots,N$. 
Then, by \eqref{proof:relation 2} and \eqref{proof:relation 1} we see in turn that 
$R_0=0$ and $R_i=0$ for $i=1,\ldots,N$. 
Therefore, we have shown that $(\eta,\mbox{\boldmath$\phi$})$ is the solution to the Isobe--Kakinuma 
model \eqref{intro:IK model}--\eqref{intro:initial conditions}.

It remains to show that the energy function $E(t)$ defined by \eqref{intro:energy} is conserved in time. 
The energy function can be written explicitly as 
\begin{align}\label{proof:energy}
E(t) &= \frac{\rho}{2}\int_{\mathbf{R}^n}\biggl\{
 \sum_{i,j=0}^N\biggl(
  \frac{1}{p_i+p_j+1}H^{p_i+p_j+1}\nabla\phi_i\cdot\nabla\phi_j
  -\frac{2p_i}{p_i+p_j}H^{p_i+p_j}\phi_i\nabla b\cdot\nabla\phi_j \\
&\makebox[6em]{}
 +\frac{p_ip_j}{p_i+p_j-1}H^{p_i+p_j-1}(1+|\nabla b|^2)\phi_i\phi_j\biggr)+g\eta^2\biggr\}{\rm d}x.
 \nonumber
\end{align}
Therefore, by the direct calculation we see that 
\begin{align*}
\frac{\rm d}{{\rm d}t}E(t)
&= \rho\biggl\{\sum_{i,j=0}^N(L_{ij}\phi_j,\partial_t\phi_i)_{L^2}
 +(g\eta+\frac12(|\mbox{\boldmath$u$}|^2+w^2),\partial_t\eta)_{L^2}\biggr\} \\
&= \rho\biggl\{\sum_{i=0}^N(H^{p_i}\partial_t\eta,\partial_t\phi_i)_{L^2}-(F_0,\partial_t\eta)_{L^2}\biggr\}
 =0.
\end{align*}
The proof of Theorem \ref{intro:theorem} is complete.

%------------------------------------------------------------------------------

\bigskip
Ryo Nemoto \par
{\sc Former Address}: \par
{\sc Department of Mathematics} \par
{\sc Faculty of Science and Technology, Keio University} \par
{\sc 3-14-1 Hiyoshi, Kohoku-ku, Yokohama, 223-8522, Japan}

\bigskip
Tatsuo Iguchi \par
{\sc Department of Mathematics} \par
{\sc Faculty of Science and Technology, Keio University} \par
{\sc 3-14-1 Hiyoshi, Kohoku-ku, Yokohama, 223-8522, Japan} \par
E-mail: iguchi@math.keio.ac.jp

\end{document}